\documentclass[11pt]{article}
\usepackage{amsmath,amssymb,amsthm,mathtools,color}
\usepackage[a4paper,margin=1in]{geometry}

\newtheorem{proposition}{\bf Proposition}[section]

\usepackage{mathtools}
\usepackage{amsfonts}
\usepackage{dsfont}
\usepackage{array}

\usepackage{dsfont} 
\usepackage{epsfig}
\usepackage{multirow}
\usepackage{float}
\usepackage{graphicx}
\usepackage{graphics}
\restylefloat{table}

\newcommand{\eins}{\mathds{1}}

\newcommand{\p}{\partial}
\newcommand{\eps}{\varepsilon}

\newcommand{\dd}{\hspace{0.0pt}{\rm d}\hspace{0.0pt}}
\newcommand{\ee}{{\rm e}\hspace{1pt}}

\newcommand{\R}{\mathbb R}
\newcommand{\N}{\mathbb N}

\newcommand{\cA}{{\mathcal A}}
\newcommand{\cB}{{\mathcal B}}
\newcommand{\cC}{{\mathcal C}}

\newcommand{\cF}{{\mathcal F}}
\newcommand{\cG}{{\mathcal G}}
\newcommand{\cH}{{\mathcal H}}

\newcommand{\cK}{{\mathcal K}}
\newcommand{\cL}{{\mathcal L}}

\newcommand{\cN}{{\mathcal N}}
\newcommand{\cO}{{\mathcal O}}

\newcommand{\cS}{{\mathfrak S}}

\newcommand{\rE}{\ensuremath{\text{E}}}
\newcommand{\rV}{\ensuremath{\text{V}}}

\newcommand{\inv}{^{-1}}
\newcommand{\trans}{^{\sf T}}

\title{A pedestrian's approach to Bayesian methods for inverse problems}

\vspace{5mm}
\author{Heiko Gimperlein, Michael Oberguggenberger\\
        {\small Universit\"{a}t Innsbruck, Engineering Mathematics}\\
         {\small Technikerstr.\ 13, 6020 Innsbruck, Austria}\\
         {\small heiko.gimperlein@uibk.ac.at, michael.oberguggenberger@uibk.ac.at}
         }

\date{}

\begin{document}

\maketitle

\noindent
{\bf Abstract:} In the past two decades, the Bayesian approach to inverse problems in infinite dimensional spaces has seen significant developments in theory, Markov chain Monte Carlo (MCMC) methods, and in applications.
The Bayesian approach in finite dimension produces posterior distributions of finite dimensional random variables. In infinite dimensions, the desired objects are random fields and their laws as measures on Hilbert or Banach spaces.
There are many textbooks and tutorials available that explain Bayesian methods and MCMC in finite dimensions. To enter the field of Bayesian inverse problems in infinite dimensions, a newcomer has to engage with a quite diverse and extensive literature. The purpose of this article is to provide a basic path to the subject
that is easily accessible to a reader with general knowledge in probability and mathematical analysis. We start with a swift introduction to the finite dimensional case. The infinite dimensional case is introduced for Gaussian priors on Hilbert spaces and Gaussian-type likelihood functions. In this situation, mere continuity of the forward map suffices to prove that the posterior measure is well-defined. We also give some new and simple stability results. Two application examples are elaborated, one finite dimensional, one infinite dimensional.
\vspace{5mm}
\vspace{5mm}\\
{\bf MSC2020 subject classifications:} Primary 65J22, 62F15; secondary 65C40, 28C20.

\tableofcontents

\section{Introduction}

Mathematical models in science and engineering are usually cast in the form of an input-output map
\[
    y = \cF(x,\theta).
\]
Here $x$ are the input variables (e.g., points in space, points of time, loads), $\theta$ are the model parameters (e.g., material or geometric parameters), $\cF$ is the input-output (or forward) map and $y$ is (a component of) the output.
If the parameters $\theta$ are known and inputs $x$ are chosen, the computation of the model response $y$
is referred to as the \emph{direct problem}.

Conversely, if the inputs $x$ and corresponding outputs $y$ are known, or given as a vector of observations $y = y_{\rm obs}$, the determination of the corresponding parameters $\theta$ is the \emph{inverse problem}.
Depending on the context, this is also referred to as model updating or parameter calibration.
The difficulty is that the inverse problem is usually ill-posed. A solution giving the required parameters $\theta$ may not exist, may be non-unique, but most importantly usually will not depend continuously on the data.
We denote the error between measured values $y_{\rm obs}$ and the theoretical model predictions $\cF(x,\theta)$ by
\[
   e = y_{\rm obs} - \cF(x,\theta).
\]
The variational approach consists in finding a set of parameters $\widehat{\theta}$ which minimizes some norm of the errors, for example, the Euclidean norm
\[
  \|e\|^2 = \|y_{\text{obs}} - \cF(x,\theta)\|^2 \rightarrow \min !
\]
In this simple form, the problem is usually unsolvable and has to be complemented by some regularization procedure such as Tikhonov regularization \cite{Engl:1996}.

In the past decades, an alternative approach that avoids the difficulties of the variational methods has gained broad acceptance: the Bayesian approach.
In a nutshell, the parameters and output variables (and possibly the input variables) are all considered as random variables $\Theta$, $Y$ and possibly $X$; the observation  $y_{\text{obs}}$ is viewed as a sample of $Y$.
Model and data uncertainties are encoded in a prior distribution $\pi_{0}(\theta)$ of $\Theta$ and the likelihood function $\pi(y_{\text{obs}}|\theta)$. The solution of the inverse problem consists in the posterior distribution of $\Theta$
\[
 \pi(\theta|y_{\text{obs}})\propto \pi(y_{\text{obs}}|\theta)\pi_{0}(\theta).
\]
Note that the Bayesian solution contains much more information on the parameter $\Theta$ than the optimization approach: its probability distribution.

In the finite dimensional case, when $x,y,\theta$ belong to finite dimensional Euclidean spaces, the Bayesian procedure, together with its numerical methods, is well-established and well-understood \cite{HeldBove:2014,Kaipio:2005,robert1999monte,Tarantola:1987}, see also the recent survey \cite{Banerjee:2026}, the
tutorial \cite{Waqar:2023} and the recent exposition \cite{Kiran:2025} from an engineering, finite element point of view.

Turning to the infinite dimensional case, let us suppress the dependence on $x$ in the notation (fixed input data may be considered as inherent to the forward map). In this situation we write the forward map as
\[
   y = \cG(\theta).
\]
The parameters $\theta$ and the output $y$ are now functions belonging to infinite dimensional (separable) Hilbert spaces $\cH,\cH'$. Given an observation $y_{\rm obs}\in \cH'$ (possibly given as a discrete version), the Bayesian inverse problem
consists in reconstructing $\theta$ as an $\cH$-valued random variable, or equivalently, constructing its law as a measure on $\cH$.

The seminal paper of Andrew Stuart \cite{Stuart_2010_acta_numerica} from 2010 sparked a tremendous development in theory and applications as well as Markov chain Monte Carlo (MCMC) methods in infinite dimensional spaces \cite{Cotter:2013}.
A key ingredient is the theory of Gaussian measures on Hilbert spaces as well as Markov chains on general state spaces.

The infinite dimensional Bayesian theory has been elaborated in more detail in various papers \cite{DashtiStuart:2017Bayesian,Lasanen:2012_1,Lasanen:2012_2}, and it even has been taken up in textbooks \cite{Nickl:2023,Sullivan:2015}.
Nevertheless, a newcomer entering the field has to engage in a quite extensive literature study. We found that easily accessible introductory expositions are lacking.

The purpose of this article -- or rather tutorial -- is to close this gap in the literature. We intend to offer an uninitiated reader
\begin{itemize}
\item a swift introduction to finite dimensional Bayesian inverse problem, together with basic MCMC algorithms;
\item a short exposition of the Bayesian approach in infinite dimensions, together with corresponding MCMC algorithms;
\item explications of the material with the help of a finite dimensional and an infinite dimensional example;
\item a comprehensive collection of results and concepts needed from Markov chains and Gaussian measures on Hilbert spaces, with proofs where useful and references to the literature, guaranteeing the existence of the posterior measure and the convergence of the MCMC algorithms.
\end{itemize}
We achieve this by restricting attention to some basic, paradigmatic situations -- only Gaussian measures on Hilbert spaces and likelihood functions with a Gaussian kernel are considered. In this case, the assumptions on the forward map $\cG$ of \cite{Stuart_2010_acta_numerica} can be dramatically
weakened; continuity suffices. Concerning the stability of the Bayesian solution under perturbations of the forward map and the observation, we obtain also new and easily accessible results, using weak convergence of measures.

Simplifications of the original assumptions in \cite{DashtiStuart:2017Bayesian,Stuart_2010_acta_numerica} have also been obtained in \cite{Latz:2023}; this encompassing paper notably contains a profound investigation of stability in various metrics on spaces of probability measures.

The plan of the paper is as follows. In Section 2, we give an introduction to the finite dimensional case. Subsection 2.1 starts with an elementary definition of conditional distributions of random variables on finite dimensional spaces. Subsection 2.2 explains the Bayesian updating procedure, Subsection 2.3 introduces
the Bayesian approach to inverse problems. Subsection 2.4 introduces the Metropolis-Hastings algorithm as an elementary Markov chain Monte Carlo method; Subsection 2.5 elaborates on some practical aspects and heuristic convergence diagnostics.
Section 3 is devoted to the generalization to infinite dimensional spaces. In Subsection 3.1, we extend the Bayesian approach to inverse problem to the Hilbert space setting with Gaussian priors and Gaussian-type likelihood functions. In this case, mere continuity of the forward map allows one to establish the existence of the posterior measure and the fact that it has finite first and second moments. The stability of the posterior measure under perturbations of the forward map and the observations is
established in Subsection 3.2, again under rather weak conditions. The posterior measures converge weakly, while the means and covariance operators converge in norm. Subsection 3.3 presents the pCN-algorithm as an infinite dimensional MCMC method.

Section 4 highlights the Bayesian methods by means of two examples. Subsection 4.1 demonstrates the reconstruction of parameters in a nonlinear spring model from observations of a finite number of maximal velocities, with visual analyses of the Markov chains. Subsection 4.2 addresses the reconstruction of the right-hand side in a two-dimensional Poisson equation from observations of the solution.

Finally, Section 5 is intended to summarize and explain a minimal amount of tools from probability theory needed for a full understanding of the procedures from Sections 2 and 3. Subsection 5.1 starts with a general introduction of conditional distributions. Section 5.2 collects basics from the theory of Markov chains and all concepts needed to infer
the existence of an invariant probability measure and to analyze convergence. Section 5.3 presents a proof of convergence of the Metropolis-Hastings algorithm in the form as applied in Subsection 2.4 as well as the ergodicity of the pCN-algorithm, based on references to the literature.
Finally, Subsection 5.4 is a short introduction to random fields and Gaussian measures on Hilbert spaces, together with the Karhunen-Lo\`eve expansion and various examples of constructing covariance operators, as applied in Subsection 4.2. We explain the probabilistic concepts in detail; however, the reader is assumed to have basic knowledge of operator theory, Sobolev spaces and weak solutions to the Poisson equation on the level of, say, Sections 1 and 2 in \cite{LordPowellShardlow:2014}.

\section{Bayesian methods: principles in finite dimensions}
\label{sec:Bayes}

This section highlights the elementary mathematical concepts needed for the Bayesian approach to inverse problems in finite dimensions. We start with an elementary version of conditional distributions and Bayesian updating, turn to Bayesian inverse problems with additive noise, and we explain one of the simplest Markov chain Monte Carlo algorithms for sampling from the posterior distribution, the Metropolis-Hastings algorithm with symmetric, everywhere positive transition kernel. Deriving the convergence properties of the algorithm requires deeper methods; an exposition of the full concepts and methods is deferred to Section \ref{Sec:appendix}.

\subsection{Bayes' formula}
\label{subsec:Bayesformula}

The notion of \emph{conditional distributions} of random variables is essential to Bayesian thinking. In the elementary, finite dimensional case only a na\"{i}ve version is needed, which can be shown to be a special case of the general notion discussed in Subsection \ref{subsec:probability}.

\emph{Discrete random variables.} Let $(\Theta,Y)$ be a bivariate discrete random variable taking values $(\theta_i,y_j)$ with probabilities $P(\Theta = \theta_i,Y = y_j)$.
The \emph{conditional distribution of $\Theta$ given $Y$} is given through the probability values
\[
  P(\Theta= \theta_i|Y = y_j) = \frac{P(\Theta = \theta_i,Y = y_j)}{P(Y=y_j)}
\]
and similarly for $P(Y = y_j|\Theta= \theta_i)$. The definitions result immediately in Bayes' formula for discrete random variables
\[
   P(\Theta= \theta_i|Y = y_j) = \frac{P(Y = y_j|\Theta= \theta_i)P(\Theta = \theta_i)}{P(Y=y_j)}.
\]

\emph{Continuous random variables.} Consider a bivariate random variable $(\Theta,Y)$ and assume that it has a probability density $\pi(\theta,y)$ with respect to Lebesgue measure $\dd\theta\times \dd y$ on $\R\times\R$ with cumulative joint distribution function $\Pi(\theta,y)$. Denote the marginal densities by
$\pi_0(\theta) = \int_{-\infty}^\infty\pi(\theta,y)\dd y$, $\pi_1(y) = \int_{-\infty}^\infty\pi(\theta,y)\dd \theta$.
Introducing the cumulative distribution function of $\Theta$ given $Y$ as
\begin{eqnarray*}
   \Pi(\theta|y) &=& \!\!\!\lim_{\Delta y \to 0}P(\Theta\leq \theta|y \leq Y \leq y+ \Delta y)
       = \lim_{\Delta y \to 0}\frac{P(\Theta\leq \theta,y \leq Y \leq y+ \Delta y)}{P(y \leq Y \leq y+ \Delta y)}\\
       & = &  \!\!\!\lim_{\Delta y \to 0} \frac{\int_{-\infty}^\theta \int_y^{y + \Delta y} \pi(\theta',y')\dd\theta'\dd y'}{\int_y^{y + \Delta y}\pi_1(y')\dd y'}
          = \frac{\int_{-\infty}^\theta \pi(\theta',y)\dd \theta'}{\pi_1(y)}
\end{eqnarray*}
and taking one more derivative with respect to $\theta$ defines the \emph{conditional density of $\Theta$ given $Y$} as
\begin{equation}\label{eq:conddens}
   \pi(\theta|y) = \frac{\pi(\theta,y)}{\pi_1(y)}
\end{equation}
provided $\pi_1(y) \neq 0$, and similarly for $\pi(y|\theta)$.
This gives rise to \emph{Bayes' formula for conditional densities}:
\begin{equation}\label{Bayesformuladensities}
\pi(\theta|y) = \frac{\pi(y|\theta)}{\pi_1(y)}\pi_0(\theta) = \frac{\pi(y|\theta)}{\int\pi(y|\theta)\pi_0(\theta)\dd\theta}\pi_0(\theta).
\end{equation}
Introducing the \emph{likelihood function} $\ell(\theta,y) = \pi(y|\theta)$, Bayes' formula is often written in short form as
\[
\pi(\theta|y) \propto \ell(\theta,y)\pi_0(\theta).
\]
Formulas \eqref{eq:conddens} and \eqref{Bayesformuladensities} generalize in the obvious way to multidimensional random variables  $(\Theta,Y)$ on $\R^d\times\R^{d'}$.

\subsection{Bayesian updating}
\label{subsec:Bayesupdating}

Let $Y$ be a random variable which depends in a way on another random variable $\Theta$. For example, $Y$ could have a probability distribution with hyperparameters $\Theta$, e.g., if $Y$ has a Gaussian distribution $\sim\cN(\mu,\sigma^2)$, $\Theta$ would be $(\mu,\sigma)$. Or $Y$ could be the output of a forward map which depends on model parameters $\Theta$. The Bayesian framework is intended to combine prior knowledge about $\Theta$ with an observation $y_{\rm obs}$ of $Y$ to improve (``update'') the information about $\Theta$ (and then also get better predictions for $Y$). The Bayesian doctrine says that the updating has to be done through Bayes' formula. Thus a further ingredient is needed: the likelihood function, the conditional distribution of $Y$ given $\Theta$. The updating consists in computing the posterior distribution of $\Theta$, the conditional distribution of $\Theta$ given $Y$. We explain the procedure in the case when $Y$ and $\Theta$ have probability densities with respect to Lebesgue measure. A general reference on Bayesian methods is \cite{HeldBove:2014}, an engineering oriented reference is \cite{Yuen:2010}.

\emph{The prior distribution.} Prior knowledge about $\Theta$ is encoded in the prior distribution with density $\pi_0(\theta)$, known from previous experiments, theoretical considerations, known statistical properties etc. If mean and variance of the prior distribution are known, $\pi_0(\theta)$ is often taken as a Gaussian. It is also common practice to use a so-called diffusive prior, a constant density on a large interval when there is little knowledge about the prior distribution or when one wants to minimize the influence of the prior distribution. Such an interval may also come from physical restrictions on the range of $\Theta$. Detailed information and advice on the choice of a prior distribution can be found in \cite[Section 6.3]{HeldBove:2014}.

\emph{The likelihood function.} The likelihood function $\ell(\theta,y)$ is meant to be proportional to the probability of observing $y$ when $\theta$ has a given value. In the case when $Y$ has a probability distribution with hyperparameters $\Theta$, the likelihood function is naturally taken as $\pi(y|\theta)$, the probability density of $Y$ with fixed parameter values $\theta$. In Bayesian inverse problems, the likelihood function is derived from the distribution of the error (see Subsection \ref{subsec:inverse}). Note that the observation $y$ could be multidimensional, e.g. a sample of $Y$ of size $n$, in which case the likelihood function is the multidimensional density of $y$ given $\theta$.

\emph{The posterior distribution.} Given an observation $y = y_{\rm obs}$, the posterior distribution is defined through its density $\pi(\theta|y) \propto \ell(\theta,y)\pi_0(\theta)$ as suggested by Bayes' formula \eqref{Bayesformuladensities}.

Schematically, the updating procedure can be described as follows. The input is

\begin{itemize}
\item a prior distribution with density $\pi_0(\theta)$;
\item a likelihood function $\ell(\theta,y)$ proportional to $\pi(y|\theta)$. The minimal requirement is that $\ell(\cdot,\cdot)$ is Lebesgue-measurable;
\item an observation $y = y_{\rm obs}$.
\end{itemize}
The updated information about $\Theta$ given the observation is contained in the posterior distribution with density
\[
   \pi(\theta|y) \propto \ell(\theta,y)\pi_0(\theta),
\]
more precisely,
\[
  \pi(\theta|y) = \frac{1}{Z(y)}\ell(\theta,y)\pi_0(\theta), \quad Z(y) = \int\ell(\theta,y)\pi_0(\theta)\dd\theta.
\]
This procedure makes sense provided $0< Z(y) < \infty$.

The factor of proportionality $Z(y)$ can be obtained by evaluating the defining integral.
If $\Theta$ is low-dimensional, the factor can be computed by numerical integration. However, we shall see in Subsection \ref{Metropolis-algorithm} that knowledge of the factor is not required for sampling from the posterior distribution.

\subsection{Bayesian inverse problems}
\label{subsec:inverse}

In this section, we show how Bayesian updating is used to solve the inverse problem for $Y = \mathcal{F}(X,\Theta)$. The crucial step is to supplement the forward model by an error term, a mean zero random variable acting as additive noise
\begin{equation}\label{eq:modelwitherror}
Y = \mathcal{F}(X,\Theta)+E.
\end{equation}
This is plausible when the error is interpreted as the measurement error of observations of $Y$. However, $E$ is often considered to subsume imprecision in the forward map $\cF$ or inaccuracies introduced by numerical approximations. The error could also be split in several summands (as is done in Subsection \ref{subsec:nonlinearspring}) or could arise in a nonlinear way (see the examples in \cite[Section 3.2]{Kaipio:2005}). For simplicity of presentation, we stay with the additive form \eqref{eq:modelwitherror}.

At this stage, we assume or postulate that the error $E$ has a known probability density $\rho(e)$. In terms of realizations $x,y,e$,  relation \eqref{eq:modelwitherror} reads $y = \cF(x,\theta) + e$. The key idea is to construct the likelihood function by substituting this relation in the error distribution:
\begin{equation}\label{eq:likelihoodfrom error}
   \ell(\theta, y) = \rho(y - \cF(x,\theta)).
\end{equation}
We apply Bayesian updating with this likelihood function. We require again
\begin{itemize}
\item a prior distribution of $\Theta$ with density $\pi_0(\theta)$;
\item an observation $y = y_{\rm obs}$.
\end{itemize}
The procedure results in the posterior distribution of $\Theta$ given $y = y_{\rm obs}$  with density
\begin{equation}\label{posterior}
    \pi(\theta|y) \propto \rho(y - \cF(x,\theta))\pi_0(\theta).
\end{equation}
From the Bayesian point of view, $\pi(\theta|y)$ constitutes the solution to the inverse problem \cite[Chapter 3]{Kaipio:2005},
as it provides the probability distribution of the parameter $\Theta$ to be reconstructed, given observation data $y_{\rm obs}$ and prior knowledge encoded in the prior density.

We specialize the procedure to the common assumption that the error has a Gaussian distribution \cite{Kaipio:2005,Tarantola:1987,Yuen:2010}. In practical applications, the observation is multidimensional, often a sample $y = (y_1,\ldots, y_n)$ obtained from model inputs $x = (x_1,\ldots,x_n)$. The error $e = (e_1,\ldots, e_n)$ then obeys an $n$-dimensional Gaussian distribution with covariance matrix $\Sigma_E$ and thus has probability density
\begin{equation}\label{error}
   {\rho(e)} \propto \exp\Big( -\frac{1}{2} e\trans \Sigma_E^{-1} e \Big).
\end{equation}
The covariance matrix $\Sigma_E$ is usually based on the empirical covariance matrix of the observed data \cite{Rosic:2016}.
If in addition the individual observations in $y_{\rm obs}$ are uncorrelated, the covariance matrix $\Sigma_{E}$ is of diagonal form $\sigma^2 I$ where $I$ is the identity matrix and $\sigma^2$ is the empirical variance of $y_{\rm obs}$, then \eqref{error} reduces to
\begin{equation}\label{error2Dreduced}
 {\rho(e)} \propto \exp\Big( -\frac{1}{2\sigma^2}\|e\|^2\Big)
\end{equation}
with the Euclidean norm $\|e\|$. Finally, the likelihood function is given by
\begin{align}\label{likelihoodGauss}
  \ell(\theta,y) = \exp\Big( -\frac{1}{2} (y-\mathcal{F}(x,\theta) )\trans \Sigma_E^{-1} (y-\mathcal{F}(x,\theta) ) \Big).
\end{align}

\emph{Bayesian estimators.} Given the posterior distribution of $\Theta$, one can obtain point estimates $\widehat{\theta}$ for $\theta$ \cite[Section 3.1.1]{Kaipio:2005}.
Commonly used estimators are the mean and the median of the posterior distribution.
The posterior mean is defined as follows:
\begin{equation}\label{eq:posterior_mean}
    \widehat{\theta} = \int \theta\, \pi(\theta|y)\,{\rm d}\theta.
\end{equation}
Given a statistical sample $\theta_1,\ldots,\theta_N$ of the posterior distribution, obtained by $N$ random and independent draws,
an approximation to \eqref{eq:posterior_mean} can be obtained as the sample mean $\frac1{N}\sum_{j=1}^N\theta_j$. Then, one can estimate the approximation error as follows.
Suppose the true value to be reconstructed is $\theta^\ast$. Then
\[
   \Big\vert\frac1{N}\sum_{j=1}^N\theta_j - \theta^\ast\Big\vert \leq \Big\vert\frac1{N}\sum_{j=1}^N\theta_j - \int \theta\, \pi(\theta|y)\,{\rm d}\theta\Big\vert
       + \left\vert\int \theta\, \pi(\theta|y)\,{\rm d}\theta - \theta^\ast\right\vert.
\]
The first summand on the right-hand side can be estimated by the usual Monte Carlo error estimate: the average error is of magnitude $\sigma/\sqrt{N}$, where $\sigma^2$ is the variance of $\Theta$. One can do slightly better asymptotically
with quasi-Monte Carlo random generators for which a pointwise error bound ${\mathcal O}(\log N/N)$ can be achieved if $\pi(\theta|y)$ is sufficiently regular \cite{Niederreiter:1992}.
{The numerical error of the numerical approximation could be included by a further splitting of the first summand, noting that the statistical sample $\theta_1,\ldots,\theta_N$ is actually replaced by a numerically approximated sample $\widetilde{\theta}_1,\ldots,\widetilde{\theta}_N$. The second summand could also be further separated into an error due to the bias of the estimator \eqref{eq:posterior_mean} and a remaining systematic error. Methods for bias correction in MCMC samples can be found in the literature \cite{Efron:2015, Jacob:2020}. There is no control on the systematic error, which subsumes errors in the choice of the likelihood function and its parameters and the choice of the physical model.}

\subsection{The Metropolis-Hastings algorithm -- Markov chain Monte Carlo}
\label{Metropolis-algorithm}

Given a probability distribution $\pi$, Markov chain Monte Carlo algorithms are designed to produce a realization $\xi_0,\xi_1,\xi_2,\ldots $ of a Markov chain which has $\pi$ as its stationary distribution.

A Markov chain is a sequence of random variables $X_0, X_1, X_2,\ldots$ on a probability space $(\Omega,\cA,P)$ with values in a state space $(S,\cS)$ with the property
-- informally speaking -- that the probability of a transition from one state to another depends only on the current state. More precisely, the conditional probability distributions (see Subsection \ref{subsec:probability}) satisfy
\[
P(X_{n+1}|X_0, X_1,\ldots, X_n) = P(X_{n+1}|X_n).
\]
A general treatment of Markov chains can be found in Subsection \ref{subsec:Markovchains}.
In this Subsection, we treat the finite dimensional case. The state space $(S,\cS)$ will be an open subset or the closure of an open subset of $\R^d$ with the $\sigma$-algebra $\cS$ of Lebesgue-measurable sets. We also assume that $\pi$ has a density with respect to Lebesgue measure, i. e. $\pi(\dd x) = \pi(x)\dd x$, and that the
finite dimensional joint distributions of the members of the chain have a probability density.
Let $p^{(n)}(y)$ be the marginal (i.e., unconditional) probability density of $X_n$.
There is a so-called transition kernel $p(x,y)$, such that at given $x$, $p(x,y)$ is a probability density and at given $y$, $p(x,y)$ is a measurable function of $x$.
The marginal probability densities are computed as
\[
   p^{(n)}(y) = \int_S p(x,y)p^{(n-1)}(x){\rm d}x
\]
and thus are uniquely determined by the initial probability distribution $p^{(0)}$.  A probability distribution $\pi$ is called invariant (or stationary) with respect to the Markov chain if
\[
\pi(y) = \int_S p(x,y)\pi(x){\rm d}x.
\]
Under rather general assumptions (see Subsection \ref{subsec:Markovchains}), a Markov chain has a unique invariant probability distribution and the distributions $p^{(n)}$ of $X_n$ converge to the stationary distribution $\pi$ in the following sense: Let $p^{(0)}$ be any initial probability distribution and $f(x)$ be any function integrable with respect to $\pi(x){\rm d}x$. Then
\[
  \lim_{N\to\infty} \frac1{N} \sum_{n=1}^N f(\xi_n) = \int_S f(x)\pi(x){\rm d}x
\]
for almost all realizations $\xi_0,\xi_1,\xi_2,\ldots $ of the chain, see Proposition \ref{prop:ET} for a precise formulation.
This is often referred to as the ergodic theorem. Accordingly, the end pieces $\xi_M, \dots, \xi_N$ for sufficiently large $M$ (after the so-called burn-in phase) can be treated as a sample of the limiting distribution $\pi$.
By the ergodic theorem, the sample mean, moments and quantiles of $\xi_M, \dots, \xi_N$ converge to the expectation value, moments and quantiles of $\pi$
as $N\to\infty$.

We present one of the most popular MCMC-methods, the \emph{Metropolis-Hastings algorithm}.
As preparation for the algorithm, one has to choose a transition kernel $q(x,y)$ on $S$, the so-called proposal distribution. We make the (standard) simplifying assumption that the kernel is symmetric and everywhere positive
\begin{equation}\label{eq:propdist}
  q(x,y) = q(y,x) > 0, \quad \text{for\ all\ } x,y\in S.
\end{equation}
Otherwise, the choice of $q$ is arbitrary.
In practice, the proposal distribution is often taken of the form $q(x,y) = q(x-y)$, where $q$ is the density of a multivariate Gaussian.

Next, one has to choose an initial distribution $p^{(0)}$ (possibly the Dirac mass at a point $x$). The Metropolis-Hastings algorithm is executed as follows.
\begin{itemize}
    \item Sample a value $\xi_0$ from $p^{(0)}$, or simply set $\xi_0 = x$ for some $x\in S$.
    \item For $k = 1,\dots, N$, sample a candidate value $\eta$ from the proposal distribution $ q(\xi_{k-1},\cdot)$.
    \item Accept/reject step. Set
    \[
      \xi_{k} = \left\{\begin{array}{ll}
          \eta & \displaystyle\text{with\ probability\ }\ a(\xi_{k-1},\eta) = \min\left(1,\frac{\pi(\eta)}{\pi(\xi_{k-1})}\right)\\[6pt]
          \xi_{k-1} & \text{with\ probability\ }\ 1 -a\big(\xi_{k-1},\eta\big).
          \end{array}\right.
    \]
\end{itemize}
The function $a(u,v) = \min\left(1,\frac{\pi(v)}{\pi(u)}\right)$ is called the \emph{acceptance probability}. If $\xi_k = \eta$ one says that the proposed candidate $\eta$ is accepted; if $\xi_k = \xi_{k-1}$ it is rejected.
The accept/reject step can be effected as follows:
\begin{itemize}
\item[] Compute $\displaystyle r = \frac{\pi(\eta)}{\pi(\xi_{k-1})}$.
\item[] Draw a random number $\zeta$ from the uniform distribution on $[0, 1]$.
\begin{itemize}
    \item If $\zeta\leq r,$ set $\xi_k=\eta$.
   \item If $\zeta > r,$ set $\xi_k=\xi_{k-1}$.
\end{itemize}
\end{itemize}
The advantage of the Metropolis-Hastings algorithm is that it requires knowledge of $\pi$ only up to multiplicative constants (only the ratio
$r$ matters). Therefore, it is particularly useful in computing the posterior distribution in a Bayesian framework.

According to the arguments collected in Subsection \ref{subsec:MCMCconvergence},
the Metropolis-Hastings chain has $\pi$ as its invariant probability distribution and satisfies the ergodic theorem.
In particular, the end pieces $\xi_M, \dots, \xi_N$ for sufficiently large $M$ provide a sample of $\pi$.

An important indicator of the quality of the Markov chain is the acceptance rate. This is the proportion of accepted values, that is, $1/N$ times the number of $k$'s for which the candidate value $\eta$ is not rejected. If the acceptance rate is close to one, the chain tends to become the Markov process generated by the proposal distribution alone. If the acceptance rate is small, the chain hardly moves and covers the state space only very slowly. Thus a medium size acceptance rate is desirable. Studies of the efficiency of the Metropolis-Hastings algorithm suggest an optimal acceptance rate of $44\%$ for univariate chains and around $23\%$ for multivariate chains \cite{Gelman:1996}.

\subsection{Practical aspects of the Metropolis-Hastings algorithm}
\label{subsec:practicalMCMC}

The Metropolis-Hastings algorithm delivers a realization $\xi_0,\xi_1,\xi_2,\ldots $ of the Markov chain under consideration. The question is if this simulated chain is of sufficient statistical quality so that it allows one to make inferences about the actual Markov chain.

\emph{Burn-in phase.} The so-called burn-in phase is defined to be the initial segment of the simulated chain, before its elements hit the main range of the posterior distribution. The length of the burn-in phase depends on the choice of the initial value (and the variance of the proposal distribution). In the application example of Subsection
\ref{subsec:nonlinearspring}, the first 1000 elements of the chain for the components $a$ and $c$ are depicted in Figure \ref{fig:burnin}. It can be observed that the length of the burn-in phase is well below 300.

\emph{Monitoring convergence.}
The quality of the simulated realization of the Markov chain can be tested by means of visual inspection (trace plots), acceptance rate and tests for stationarity.
Additional methods of analysing the quality of Markov chains produced by MCMC can be found in
\cite[Chapter 12]{robert1999monte} and in \cite{Brooks:1998}. For further  statistical explorations of posterior distributions see, e.g., \cite{Efron:2015}
and \cite{Barndorff-Nielsen:1994}.

To assess the quality of the simulated chain, one usually checks if it explores the whole state space (the range of the parameter under consideration) reasonably fast. One important criterion is the acceptance rate, which optimally should be around 44\% for univariate chains and around $23\%$ for multivariate chains, see Subsection \ref{Metropolis-algorithm}. In the application example of Subsection \ref{subsec:nonlinearspring}, the acceptance rate turned out to be 24\%, in the infinite dimensional example of Subsection \ref{subsec:2dPoisson} it was 15\%, both reasonably close to the optimal ones.

A second criterion is whether the members of the chain move between different parts of the state space quickly. This can be checked visually and is evident for the chains in
Section \ref{sec:Applications} from Figures \ref{fig:burnin}, \ref{fig:chains_histograms} and \ref{fig:pyramid_chains_histograms}.
To check this and other properties for multi-dimensional chains, a host of graphical and non-graphical methods has been proposed. We refer to the survey articles \cite{Cowles:1996,Roy2020}.

The final question concerns the convergence of the chain to its stationary distribution. In terms of the given realization $\xi_0,\xi_1,\xi_2,\ldots $, this means that subsamples $\xi_K, \ldots, \xi_{K+M}$ and $\xi_L, \ldots, \xi_{L+M}$, with sufficiently large $K$ and $L>K+M$, should have the same distribution. This can be tested by nonparametric statistical tests such as the Kolmogorov-Smirnov test. However, these tests apply only to iid (independent and identically distributed) subsamples.  As suggested in \cite[Section 12.2.2]{robert1999monte}, independent sampling can be achieved by selecting batches of the form
$\xi_K, \xi_{K + G}, \xi_{K+2G},\ldots, \xi_{K+lG}$ with $G >\!\!> 1$ (and $lG\leq M$).

\emph{Confidence intervals.}
Approximate confidence intervals for the mean can be computed using the central limit theorem for Markov chains (Proposition \ref{prop:CL}) from Subsection \ref{subsec:Markovchains}.
Let $S_N = X_1 + \cdots + X_N$. If the random variables $X_j$ were iid with variance $\sigma^2$, the variance of $S_N$ would be $N\sigma^2$. When the $X_j$ arise as members of a Markov chain,
the variance of $S_N$ can still be asymptotically estimated.
Under certain conditions, ${\rm V}(S_N)/N$ converges to the number $\gamma^2$ given by \eqref{eq:gamma2} as $N\to\infty$.
Further, if $\gamma > 0$, the random variable $S_N/\gamma\sqrt{N}$ converges (in distribution) to the standard normal distribution. The required mixing conditions
(examples of such conditions are in Subsection \ref{subsec:Markovchains}) cannot easily be verified in practice, but can be enforced in the MCMC setting (see the discussion in \cite[Sections 6.7.2 and 6.9]{robert1999monte}).

Thus, taking a segment $\xi_{M+1}, \ldots, \xi_{M+N}$ of length $N$ of the chain after the burn-in period with sample mean $\mu_{N}= \frac{1}{N}\sum_{i=1}^N\xi_{M+i}$, one may assume that
$\mu_{N}$ is approximately normally distributed around $\mu^\ast$ with variance $\gamma^2/N$; here $\mu^\ast$ is the true mean of the posterior distribution.
Due to the assumed stationarity of the chain, the variance and the covariances in \eqref{eq:gamma2} can be estimated by the empirical autocovariance function of the segment $\xi_{M+1}, \ldots, \xi_{M+N}$, and the infinite sum can be truncated at large $k = K$ when the covariances become negligible (empirically around $K = 10$ to $20$).

Based on these considerations, approximate confidence intervals for the posterior mean can be computed as $\mu_{N} \pm q\gamma/\sqrt{N}$, where $q$ is the quantile of the standard normal distribution corresponding to the desired confidence level ($q=1.96$ for 95\% confidence).

This procedure has been applied, for example, in \cite[Figure 5]{GMMO} to observe the behavior of confidence intervals with increasing $N$.

\section{Bayesian methods: infinite dimensions}
\label{sec:BayesInf}

This section requires some basics of the theory of Gaussian measures on Hilbert spaces, as collected in Subsection \ref{subsec:Gaussianmeasures}.

\subsection{Bayesian inverse problems in Hilbert spaces}
\label{subsec:inverseHilbert}

The material in this section is based on the seminal paper \cite{Stuart_2010_acta_numerica} and on \cite{DashtiStuart:2017Bayesian}. However, we do not need the elaborate general hypotheses as in \cite[Assumptions 2.6 and 2.7]{Stuart_2010_acta_numerica}, because we consider only likelihood functions of the form \eqref{eq:likelihoodinfdim} below. For the first step in constructing the posterior measure, only continuity of the forward map $\cG$ is required; for convergence of approximations, polynomial growth conditions are needed, in addition.
We do not need local Lipschitz continuity of the forward map, as often required \cite{DashtiStuart:2017Bayesian,Nickl:2023,Stuart_2010_acta_numerica,Sullivan:2015}. Similar conditions on the forward map without Lipschitz continuity have also been stated in \cite{Latz:2023}.

Let $\cH,\cH'$ be separable Hilbert spaces and consider the forward problem $Y = \cG(\Theta)$ for a map $\cG:\cH\to\cH'$. The inverse problem consists in reconstructing the random variable $\Theta$ from an observation $y_{\rm obs}$ of $Y$. Its solution follows the procedure of Bayesian updating described in Subsection \ref{subsec:Bayesupdating}. We require the following
\begin{itemize}
\item a prior probability measure $\pi_0(\dd\theta)$ on $\cH$. The prior is assumed to be a nondegenerate Gaussian measure with  covariance operator $\cC_0$,
${\rm Ker}(\cC_0) = \{0\}$;
\item a likelihood function $\ell(\theta,y)$ on $\cH\times\cH'$. We assume that the likelihood function is of the form
\begin{equation}\label{eq:likelihoodinfdim}
\ell(\theta,y) = \exp\Big(-\frac1{2\sigma^2}\|y - \cG(\theta)\|_{\cH'}^2\Big)
\end{equation}
with $\sigma^2$ interpreted as variance of the error.
\item an observation $y = y_{\rm obs}\in\cH'$.
\end{itemize}
In analogy with the finite dimensional case, the solution of the Bayesian inverse problem is the probability measure $\pi^y(\dd\theta)$ defined by
\begin{equation}\label{eq:solutionIP}
\pi^y(\dd\theta) = \frac{1}{Z(y)}\ell(\theta,y)\pi_0(\dd\theta), \quad Z(y) = \int\ell(\theta,y)\pi_0(\dd\theta).
\end{equation}
In some cases, it can be shown that $\pi^y(\dd\theta)$ is indeed the conditional distribution of $\Theta$ given $y$, see Remark (a) below as well as the general discussion in \cite{Lasanen:2012_1}.

To prove that this is a valid procedure, we have to show that $\pi^y(\dd\theta)$ is a well-defined probability measure on $\cH$, that is, $0<Z(y)<\infty$.

\begin{proposition}\label{prop:postwelldefined}
Assume that $\cG:\cH\to\cH'$ is continuous and let $y\in\cH'$. Then\\
(a) The measure $\pi^y(\dd\theta)$ given by \eqref{eq:solutionIP} is a well-defined probability measure.\\
(b) The measure $\pi^y(\dd\theta)$ has a well-defined mean  and covariance operator
\[
\langle m, f\rangle = \int\langle \theta, f\rangle \frac{1}{Z(y)}\ell(\theta,y)\pi_0(\dd\theta),\quad
\langle \cC f, g\rangle = \int\langle \theta-m, f\rangle \langle \theta-m, g\rangle \frac{1}{Z(y)}\ell(\theta,y)\pi_0(\dd\theta)
\]
for $f,g\in\cH$.
\end{proposition}
\emph{Proof.} (a) First, the continuity of $\cG$ implies that  $\ell(\theta,y)$ is measurable with respect to $\theta$.
Second, $\|y - \cG(\theta)\|_{\cH'}^2 \geq 0$, so $\ell(\theta,y) \leq 1$ and $Z(y) \leq \int\pi_0(\dd\theta) = 1 <\infty$, yielding the upper bound.
Third, the lower bound follows from the continuity of $\cG$ at zero, according to which there is $r>0$ such that $\|\cG(\theta)\|_{\cH'} \leq 1 + \|\cG(0)\|_{\cH'}$
whenever $\|\theta\|_\cH < r$. In particular, $\|y - \cG(\theta)\|_{\cH'}$ is bounded for fixed $y\in\cH'$ and so $\ell(\theta,y) \geq \exp(-K)$ for $\|\theta\|_\cH < r$ and some $K>0$. We conclude that
\[
   Z(y) \geq \int_{\{\|\theta\|_\cH < r\}} \exp(-K) \pi_0(\dd\theta) = \exp(-K)\pi_0(\{\|\theta\|_\cH < r\}) > 0,
\]
because any nondegenerate Gaussian measure is full, i. e., $\pi_0(B) > 0$ for any non-empty open subset of $\cH$ \cite[Proposition 1.25]{DaPrato:2006}.

(b) Since $\pi_0$ is a Gaussian measure, $\int \|\theta\|_\cH^p\pi_0(\dd\theta)$ is finite for every $p\in\N$ \cite[Theorem 6.6]{Stuart_2010_acta_numerica}. Since $|\langle \theta, f \rangle| \leq \|\theta\|\|f\|$, the integral defining $\langle m, f\rangle$ is finite. A similar argument applies to the quadratic form defining the covariance. \qed

\emph{Remarks.} (a) Assume that $\cH'$ is finite dimensional and $\ell(\theta,y)$ in \eqref{eq:likelihoodinfdim} \emph{is} the conditional density of $Y$ given $\Theta$. Then it can be shown that $\pi^y(\dd\theta)$ is indeed the conditional distribution of $\Theta$ given $Y$ \cite[Theorem 6.31]{Stuart_2010_acta_numerica}.

(b) The direct generalization of formula \eqref{likelihoodGauss} to the Hilbert space setting would be to let
\[
\ell(\theta,y) = \exp\Big(-\frac12\left\langle\cK^{-1/2}\big(y - \cG(\theta)\big),\cK^{-1/2}\big(y - \cG(\theta)\big)\right\rangle\Big)
\]
for some covariance operator $\cK$ on $\cH'$. However, this makes sense only when $y-\cG(\theta)$ belongs to the domain of $\cK^{-1/2}$, which is the range of $\cK^{1/2}$.
But the latter range is the Cameron-Martin space of the measure $\mu\sim\cN(0,\cK)$, which is known to have $\mu$-measure zero \cite[Lemma 6.10]{Stuart_2010_acta_numerica},
\cite[Proposition 1.27]{DaPrato:2006}. For this reason, we consider only likelihood functions of the special form \eqref{eq:likelihoodinfdim}.

\subsection{Stability of the Bayesian solution}
\label{subsec:StabilityBayes}

In this subsection, we address the question of stability of the method with respect to perturbations $\cG_n$ of the forward map $\cG$ or the observation. It will be seen that under additional conditions on $\cG_n$ and $\cG$, the corresponding probability measures $\pi^{y_n}_n(\dd\theta)$ converge to
$\pi^y(\dd\theta)$ weakly in the sense of measures (equivalently, in the Prohorov metric \cite[Theorem 6.8]{Billingsley:1999}). Further, the means and the covariance operators converge in norm, and rates of convergence can be obtained. Concentrating on convergence in the weak sense of measures allows us to considerably weaken the assumptions stated in \cite{Stuart_2010_acta_numerica}.
As mentioned in the Introduction, a thorough investigation of stability in various metrics and the corresponding hypotheses is in \cite{Latz:2023}.

Suppose a sequence of observations $y_n\to y$ in $\cH'$ and a sequence of continuous forward maps $\cG_n$ converging pointwise to $\cG$, $\cG_n(\theta)\to \cG(\theta)$ for every $\theta\in \cH$ are given. We obtain a sequence of likelihood functions
$$
\ell_n(\theta,y_n) = \exp\Big(-\frac1{2\sigma^2}\|y_n - \cG_n(\theta)\|_{\cH'}^2\Big)
$$
according to \eqref{eq:likelihoodinfdim} and a sequence of probability measures
\begin{equation}\label{eq:approx}
\pi^{y_n}_n(\dd\theta) = \frac{1}{Z_n(y_n)}\ell_n(\theta,y_n)\pi_0(\dd\theta), \quad Z_n(y_n) = \int\ell_n(\theta,y_n)\pi_0(\dd\theta).
\end{equation}
The question is in what sense the measures $\pi^{y_n}_n(\dd\theta)$ converge to  $\pi^y(\dd\theta)$ given by \eqref{eq:solutionIP}.
For clarity, we collect the hypotheses again:
\begin{itemize}
\item[(A1)] $\cG_n, \cG : \cH\to\cH'$ are continuous maps such that $\cG_n(\theta)\to \cG(\theta)$ in $\cH'$ as $n\to\infty$ for every $\theta \in \cH$; $y_n\to y$ in $\cH'$ as $n\to\infty$.
\end{itemize}
We make the following assumption.
\begin{itemize}
\item[(A2)] There are $r>0$ and $K > 0$ such that $\|\cG(\theta)\|_{\cH'} \leq K$, $\|\cG_n(\theta)\|_{\cH'} \leq K$ for all $n\in\N$ and $\|\theta\|_\cH < r$.
\end{itemize}
\begin{proposition}\label{prop:conv}
(a) Under assumptions (A1) and (A2), the posterior measures $\pi^{y_n}_n(\dd\theta)$ given by \eqref{eq:approx} converge to $\pi^y(\dd\theta)$ weakly in the sense of measures.\\[1pt]
(b) The means $m_n$ converge to $m$ in norm and the covariance operators converge in operator norm.
\end{proposition}
\emph{Proof.}
(a) Assumption (A2) and the fact that the sequence $y_n$ is bounded implies that there is a constant $K'$ such that
\[
    \ell(\theta,y) \geq  \exp(-K'), \quad \ell_n(\theta,y_n) \geq  \exp(-K')
\]
for all $n\in\N$, whenever $\|\theta\|_\cH < r$.
As in the proof of Proposition \ref{prop:postwelldefined}, this gives uniform bounds from below for $Z(y)\geq a > 0$ and all $Z_n(y_n)\geq a > 0$.

Now let $h:\cH\to \R$ be a continuous function such that $|h(\theta)| \leq c(1+\|\theta\|_\cH^p)$ for some $c>0$, $p \geq 0$ and all $\theta\in\cH$. First note that
\[
  Z_n(y_n) = \int\exp\Big(-\frac1{2\sigma^2}\|y_n - \cG_n(\theta)\|_{\cH'}^2\Big)\pi_0(\dd\theta)
    \to \int\exp\Big(-\frac1{2\sigma^2}\|y - \cG(\theta)\|_{\cH'}^2\Big)\pi_0(\dd\theta) = Z(y)
\]
as $n\to\infty$ by Lebesgue's dominated convergence theorem, since the integrands converge pointwise and are bounded by 1. Similarly,
\begin{equation}\label{eq:convergence}
\begin{array}{rcl}
   \displaystyle\int h(\theta)\pi^{y_n}_n(\dd\theta) & = & \displaystyle \int h(\theta)\frac{1}{Z_n(y_n)}
        \exp\Big(-\frac1{2\sigma^2}\|y_n - \cG_n(\theta)\|_{\cH'}^2\Big)\pi_0(\dd\theta)\\
        &&\\
   & \to & \displaystyle\int h(\theta)\frac{1}{Z(y)}\exp\Big(-\frac1{2\sigma^2}\|y - \cG(\theta)\|_{\cH'}^2\Big)\pi_0(\dd\theta) \ =\  \int h(\theta)\pi^y(\dd\theta)
    \end{array}
\end{equation}
because the integrands converge pointwise and are bounded by the $\pi_0$-integrable function $c(1+\|\theta\|_\cH^p)/a$ (proof of part (b) in Proposition \ref{prop:postwelldefined}).

In particular, \eqref{eq:convergence} holds for any bounded continuous function $\cH\to\R$, implying the weak convergence of $\pi^{y_n}_n(\dd\theta)$ to $\pi^y(\dd\theta)$.
Taking $h(\theta) = \langle\theta, f\rangle$ for $f\in\cH$, we see that the means converge in the sense that
\[
   \langle m_n,f\rangle = \int\langle\theta,f\rangle \pi^{y_n}_n(\dd\theta) \to \int\langle\theta,f\rangle \pi^y(\dd\theta) = \langle m,f\rangle
\]
for every $f\in\cH$. Similarly, the covariance operators converge when evaluated at fixed $f,g\in\cH$.

(b) The norm convergence of the means $m_n$ to $m$ follows from the next argument. Let $\eps > 0$. Since $\|m_n-m\|_\cH = \sup_{\|f\|_\cH \leq 1}|\langle m_n -m,f\rangle|$, we can find $f_n$, $\|f_n\|_\cH \leq 1$, such that
\[
   \|m_n-m\|_\cH - \frac{\eps}2 \leq |\langle m_n -m,f_n\rangle|.
\]
But $|\langle m_n -m,f_n\rangle|$ equals
\begin{eqnarray*}
 \Big\vert\int \langle\theta,f_n\rangle\left(\frac{1}{Z_n(y_n)}\exp\Big(-\frac1{2\sigma^2}\|y_n - \cG_n(\theta)\|_{\cH'}^2\Big)
     - \frac{1}{Z(y)}\exp\Big(-\frac1{2\sigma^2}\|y - \cG(\theta)\|_{\cH'}^2\Big)\right)\pi_0(\dd\theta)\Big\vert\\
   \leq \int \|\theta\|_\cH\left|\frac{1}{Z_n(y_n)}\exp\Big(-\frac1{2\sigma^2}\|y_n - \cG_n(\theta)\|_{\cH'}^2\Big)
     - \frac{1}{Z(y)}\exp\Big(-\frac1{2\sigma^2}\|y - \cG(\theta)\|_{\cH'}^2\Big)\right|\pi_0(\dd\theta)
\end{eqnarray*}
As above, dominated convergence gives that the right-hand side of the inequality converges to zero. Thus, for sufficiently large $n$,
$\|m_n-m\|_\cH - \frac{\eps}2 \leq \frac{\eps}2$ and so $\|m_n-m\|_\cH \leq \eps$ as desired. In a similar way one can show that the covariance operators converge in operator norm.
\qed

Under an additional assumption on the growth of $\cG$ we can derive \emph{rates of convergence}. In addition to (A2) we make the assumption
\begin{itemize}
\item[(A3)] There is $C>0$, $q>0$ and $\psi:\N\to \R_+$ such that, for all $\theta\in\cH$ and $n\in\N$,
\[
   \|\cG(\theta)\|_{\cH'} \leq C\big(1 + \|\theta\|_\cH^q\big), \quad \|\cG_n(\theta)\|_{\cH'} \leq C\big(1 + \|\theta\|_\cH^q\big)
\]
and
\[
   \|\cG_n(\theta) - \cG(\theta)\|_{\cH'} \leq C\psi(n)\|\theta\|_\cH^q.
\]
\end{itemize}
(Of course we usually have that the rate of model approximation $\psi(n)$ tends to zero as $n\to\infty$; compare also with the more general version
\cite[Theorem 4.6]{Stuart_2010_acta_numerica}.) Assumption (A3) allows us to prove a rate-dependent version of \eqref{eq:convergence}.

\begin{proposition}\label{prop:rates}
Assume (A1) -- (A3) and let $h:\cH\to\R$ be measurable with respect to the Borel $\sigma$-algebras on $\cH$ and $\R$ such that $|h(\theta)| \leq c(1+\|\theta\|_\cH^p)$ for some $c>0, p \geq 0$ and all $\theta\in\cH$. Then there is a constant $C''$ such that
\[
\Big\vert\int h(\theta)\pi^{y_n}_n(\dd\theta) - \int h(\theta)\pi^y(\dd\theta)\Big\vert
    \leq  C''\big(\| y - y_n\|_{\cH'} + \psi(n)\big).
\]
\end{proposition}
\emph{Proof.}
Making use of the inequality
\[
   |\ee^{-b} - \ee^{-a}| \leq |b-a|\quad\text{for}\ a,b\geq 0,
\]
we can estimate
\begin{eqnarray*}
|\ell_n(\theta,y_n) - \ell(\theta,y)| & = &
                          \Big\vert\exp\Big(-\frac1{2\sigma^2}\|y_n - \cG_n(\theta)\|_{\cH'}^2\Big) - \exp\Big(-\frac1{2\sigma^2}\|y - \cG(\theta)\|_{\cH'}^2\Big)\Big\vert\\
                         & \leq & \frac1{2\sigma^2}\Big\vert \|y_n - \cG_n(\theta)\|_{\cH'}^2 - \|y - \cG(\theta)\|_{\cH'}^2\Big\vert\\
       & = & \frac1{2\sigma^2}\Big\vert \|y_n - \cG_n(\theta)\|_{\cH'} - \|y - \cG(\theta)\|_{\cH'}\Big\vert\Big(\|y_n - \cG_n(\theta)\|_{\cH'} + \|y - \cG(\theta)\|_{\cH'}\Big)
\end{eqnarray*}
Assuming that $\|y_n\|_{\cH'}, \|y\|_{\cH'}\leq r$, we can bound the second factor on the right-hand side by
\[
\|y_n - \cG_n(\theta)\|_{\cH'} + \|y - \cG(\theta)\|_{\cH'} \leq 2\Big(r + C\big(1 + \|\theta\|_\cH^q\big)\Big).
\]
Using the reverse triangle inequality, the first factor is estimated by
\begin{eqnarray*}
\big\vert \|y_n - \cG_n(\theta)\|_{\cH'} - \|y - \cG(\theta)\|_{\cH'}\big\vert &\leq& \|y_n - \cG_n(\theta) - (y - \cG(\theta))\|_{\cH'} \\
&\leq & \| y - y_n\|_{\cH'} + \|\cG_n(\theta) - \cG(\theta)\|_{\cH'} \\
& \leq & \| y - y_n\|_{\cH'} + C\psi(n)\|\theta\|_\cH^q.
\end{eqnarray*}
Since all powers of $\|\theta\|_{\cH}$ are integrable with respect to $\pi_0(\dd\theta)$, we get
\[
  |Z_n(y_n) - Z(y)| \leq \int\big\vert\ell_n(\theta,y_n) - \ell(\theta,y)\big|\pi_0(\dd\theta)   \leq C'\big(\| y - y_n\|_{\cH'} + \psi(n)\big)
\]
for some constant $C'$. Since $Z_n(y_n)$ and $Z(y)$ are bounded from below by some constant $a>0$, we also have
\[
\Big\vert\frac1{Z_n(y_n)} - \frac1{Z(y)}\Big\vert \leq \frac{C'}{a^2}\big(\| y - y_n\|_{\cH'} + \psi(n)\big).
\]
Given a Borel measurable function $h:\cH\to\R$ such that $|h(\theta)| \leq c(1+\|\theta\|_\cH^p)$, we combine the estimates to get
\begin{eqnarray*}
\Big\vert\int h(\theta)\pi^{y_n}_n(\dd\theta)\!\!\! &-&\!\!\! \int h(\theta)\pi^y(\dd\theta)\Big\vert
    \ = \ \Big\vert\int h(\theta)\left(\frac{1}{Z_n(y_n)}\ell_n(\theta,y_n) - \frac1{Z(y)}\ell(\theta,y)\right)\pi_0(\dd\theta)\Big\vert\\[4pt]
    & \leq & \int |h(\theta)|\left(\frac{1}{Z_n(y_n)}\big\vert\ell_n(\theta,y_n) - \ell(\theta,y)\big\vert + \Big\vert\frac1{Z_n(y_n)} - \frac1{Z(y)}\Big\vert\big\vert\ell(\theta,y)\big\vert\right)\pi_0(\dd\theta)\\[4pt]
    & \leq &  C''\Big(\| y - y_n\|_{\cH'} + \psi(n)\Big)
\end{eqnarray*}
where we have used the previous estimates, the lower bound $a$ for $Z_n(y_n)$ and the fact that $\ell(\theta, y) \leq 1$.
This is the desired rate dependent version of \eqref{eq:convergence}. \qed

In numerical approximations of $\pi^y(\dd\theta)$, an additional error is introduced by truncating the Karhunen-Lo\`eve expansion to a finite number of terms.
In this respect, convergence results and error estimates can be found in \cite[Section 2.3]{DashtiStuart:2011}.
For further convergence results from different perspectives we refer to \cite{Lasanen:2012_2,Nickl:2023}.

Proposition \ref{prop:rates} is especially relevant for finite element approximations $\cG_n$ of the solution operator $\cG$ to elliptic partial differential equations, where one has convergence rates in the appropriate Sobolev norms, see e.g. \cite[Chapter 2]{LordPowellShardlow:2014}.

\subsection{MCMC for measures}
\label{subsec:MCMCfunctions}

The Metropolis-Hastings algorithm presented in Subsection \ref{Metropolis-algorithm} was designed to produce a realization $\xi_0,\xi_1,\xi_2,\ldots $ of a Markov chain on a state space $S$, assumed to be an open subset or the closure of an open subset of $\R^d$, such that its end pieces form a sample of a given probability distribution defined by a density $\pi(x)$ (with respect to Lebesgue measure).

Here the state space is a Hilbert space $\cH$, and the sought after probability distribution is a measure $\pi$ which is absolutely continuous with respect to a given Gaussian measure $\pi_0$ on $\cH$,
with Radon-Nikodym derivative proportional to some $\ell(\theta) > 0$,
\[
  \frac{\dd \pi}{\dd \pi_0}(\theta) \propto \ell(\theta).
\]
The Markov chain will deliver a sample of an $\cH$-valued random variable $\Theta$ whose law is $\pi$. We assume that $\pi_0 \sim \cN(0,\cC_0)$ with covariance operator $\cC_0$.

Main references for this section are \cite{Cotter:2013,DashtiStuart:2017Bayesian}. We present here just one of the MCMC-algorithms, the preconditioned Crank-Nicolson (pCN) algorithm.
\begin{itemize}
    \item Pick a starting value $\xi_0\in\cH$.
    \item For $k = 1,\dots, N$, propose a candidate $\eta = (1-\beta^2)^{1/2}\xi_{k-1} + \beta w_k$, where $w_k\sim \cN(0,\cC_0)$.
     \item Accept/reject step. Set
    \[
      \xi_{k} = \left\{\begin{array}{ll}
          \eta & \displaystyle\text{with\ probability\ }\ a(\xi_{k-1},\eta) = \min\left(1,\frac{\ell(\eta)}{\ell(\xi_{k-1})}\right)\\[6pt]
          \xi_{k-1} & \text{with\ probability\ }\ 1 -a\big(\xi_{k-1},\eta\big).
          \end{array}\right.
    \]
\end{itemize}

The accept/reject step is effected in the same way as in the finite dimensional case (Subsection \ref{Metropolis-algorithm}). Here $0 < \beta < 1$ is an adjustable parameter, which affects the behavior of the chain (as can be seen by considering the extreme cases $\beta =0$, $\beta= 1$).

\emph{Remarks.} (a) A realization of $w_k$ is simply obtained by drawing from the Karhunen-Lo\`eve expansion of $\cC_0$. Other sampling methods are described in Example 4 and  Remark (d) at the end of Subsection \ref{subsec:Gaussianmeasures}.

(b) The relation to the finite dimensional algorithm can be illustrated as follows. Suppose the transition kernel $q(x,y)$ is of the form $q(x,y) = q(x-y)$ where $q(\cdot)$ is the density of a centered Gaussian variable with covariance matrix $C$. Then the candidate $\eta$ in Subsection \ref{Metropolis-algorithm} is simply obtained by setting
$\eta = \xi_{k-1} + w_k$, where $w_k$ is a draw from $\cN(0,C)$. This is usually referred to as the \emph{standard random walk proposal}. It is argued in
\cite[Sections 4.2 and 6.1]{Cotter:2013} that this proposal is not suitable in infinite dimension; rather the modified version with the factor $\beta$ should be used to obtain the required ergodic properties of the Markov chain.

(c) The pCN-method is derived from a \emph{Langevin algorithm}, starting from a Langevin-type diffusion equation having $\pi_0$ or $\pi$ as invariant measure, discretizing the equation and then inserting a Metropolis-Hastings accept/reject step to produce the desired Markov chain from the discrete time step solution. In the infinite dimensional case, two such Langevin equations were considered in \cite{Beskos:2008}, one in plain and one in preconditioned form. Discretization by an implicit Euler scheme led the authors to the algorithm PIA (preconditioned implicit algorithm). Motivated by the fact that the equation contains the covariance operator and its inverse (commonly a second order elliptic differential operator), the method was later renamed to pCN in \cite{Cotter:2013}, see also \cite{HuYaoLi:2017} for a brief discussion.

(d) In the literature, the Radon-Nikodym derivative is usually written in terms of a potential $\Phi(\theta)\geq 0$ as $\ell(\theta) = \exp(-\Phi(\theta))$. One may reformulate the accept-reject step by saying that the acceptance probability for the candidate $\eta$ is $a(\xi_{k-1},\eta)$ where
\begin{equation}\label{eq:acceptanceprob}
   a(u,v) = \min\big(1,\exp(\Phi(u) - \Phi(v))\big) \quad \big( = \min(1,\ell(v)/\ell(u))\ \big).
\end{equation}

The question of convergence of the pCN-algorithm is addressed in Subsection \ref{subsec:MCMCconvergence}. It will be seen that the pCN-chain satisfies the detailed balance condition and has certain ergodicity properties.


\section{Applications}
\label{sec:Applications}

\subsection{Reconstruction of parameters in a nonlinear spring model}
\label{subsec:nonlinearspring}

The first example illustrates finite dimensional Bayesian inversion in a nonlinear model. We consider the nonlinear spring
\begin{equation}\label{eq:NLspring}
\begin{array}{l}
 \ddot{x}(t) + 10\,x(t) + c\, x(t)|x(t)|^{a-1} = 0,\ t\in[0,T],\\
  x(0) = x_0,\ \dot{x}(0) = \dot{x}_0.
 \end{array}
\end{equation}
The value of the spring constant $k=10$ is chosen for convenience. The goal is to reconstruct the coefficient $c$ in front of the nonlinear term and the exponent $a$.
To set up the forward (and observational) model, recall that the energy
\[
   E = \frac12\dot{x}^2 + 5x^2 + \frac{c}{a+1}|x|^{a+1}
\]
is conserved. From the phase diagram one reads off that the maximal velocity $\dot{x}_{\text{max}}$ occurs when $x=0$. Thus $\dot{x}_{\text{max}} = \sqrt{2E}$.
Specializing to the case when the spring is at rest initially ($\dot{x}_0 = 0$) and combining this with energy conservation gives the relation
\[
   \dot{x}_{\text{max}} = \sqrt{2}\left(5x_0^2 + \frac{c}{a+1}|x_0|^{a+1}\right)^{1/2}.
\]
Thus at fixed parameter values $\theta = (a,c)$, $Y = \dot{x}_{\text{max}}$ is a function of $X = x_0$. This will serve as our forward model $Y = \mathcal{F}(X,\Theta)$.
The forward map is evaluated by computing a trajectory $x(t)$ with given $x_0$ and calculating $\dot{x}_{\text{max}} = \max_{t\in[0,T]} \dot{x}(t)$.
The trajectories are computed by a modified forward Euler scheme applied to the $(2\times 2)$-system for $(x,\dot{x})$ where each component is successively updated, formally
\begin{equation}\label{eq:Euler}
\begin{array}{rcl}
x_{i+1} &=& x_i + \Delta t\, \dot{x}_i,\\
\dot{x}_{i+1} &=& \dot{x}_i - \Delta t\,\big(x_{i+1} + c\,x_{i+1}|x_{i+1}|^{a-1}\big).
\end{array}
\end{equation}

\begin{figure}[h]
    \centering
    \includegraphics[width = 5.5cm]{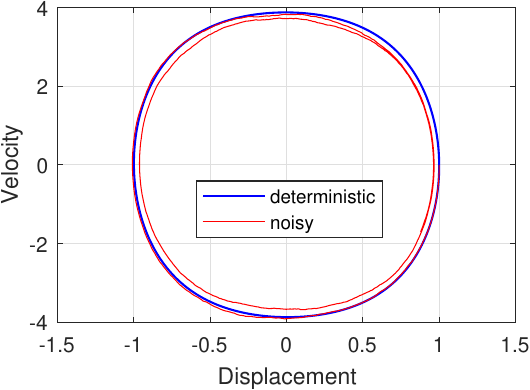} \qquad
    \caption{True and perturbed trajectory for $x_0 = 1$.}
    \label{fig:trajectory}
\end{figure}

\emph{Observed data.} The observed data will be artificially generated as stochastically perturbed solutions with reference parameters $a_\ast = 3, c_\ast = 10$, that is, as solutions to the It\^o stochastic differential equation
\begin{equation}\label{eq:NLspringIto}
\begin{array}{l}
 \ddot{x}(t) + 10\,x(t) + 10\, x(t)|x(t)|^2 = \sigma_{\text{noise}}\,\dot{B}_t,\ t\in[0,T],\\
  x(0) = x_0,\ \dot{x}(0) = 0,
 \end{array}
\end{equation}
where $B_t$ is standard Brownian motion. The stochastic differential equation is solved numerically by the Euler-Maruyama scheme corresponding to \eqref{eq:Euler}.
In the present example, we took $T = \pi$, $\Delta t = \pi/1000$, $\sigma_{\text{noise}} = 0.1$ and the seven values of $x_0 = 0.5, 1.0, 1.5, \ldots, 3.5$. The observed data $y_{\text{obs}}$ consisted of the corresponding seven maximal velocities obtained with seven random realizations of Brownian motion.
Typical trajectories of the unperturbed system \eqref{eq:NLspring} and the noisy system \eqref{eq:NLspringIto} with $x(0) = 1$, $\dot{x}(0) = 0$ are shown in Figure \ref{fig:trajectory}.

\begin{figure}[h]
    \centering
    \includegraphics[width = 5cm]{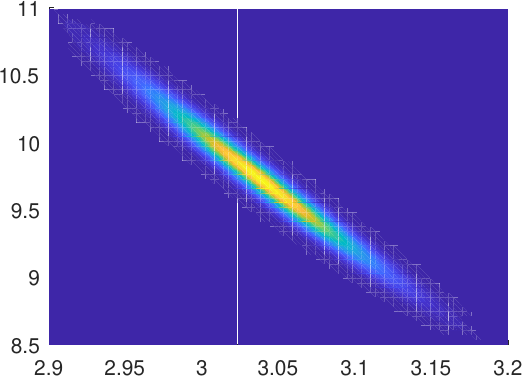}\qquad
    \raisebox{-4mm}{\includegraphics[width = 6.2cm]{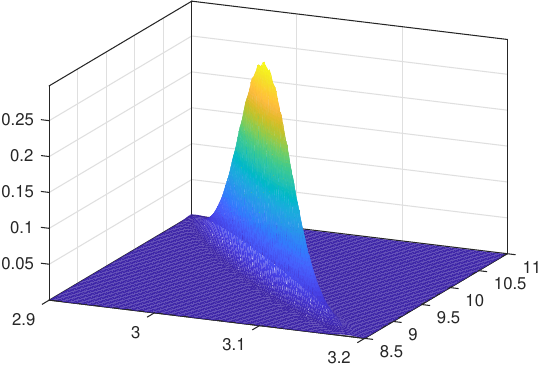}}
    \caption{Contour plot and 3D-plot of likelihood function.}
    \label{fig:likelihood}
\end{figure}

\emph{The error distribution.} We allowed for two types of error -- a model error $E_1$ due to the randomness of the trajectories of \eqref{eq:NLspringIto} and a measurement error $E_2$ in the values of the maximal velocities. We considered the errors as additive, independent and Gaussian distributed, resulting in the model
\[
   Y = \mathcal{F}(X,\Theta) + E_1 + E_2.
\]
The covariance matrix $\Sigma_{E_1}$ was estimated by running the model \eqref{eq:NLspringIto} 500 times and taking the empirical covariance matrix of the 500 observation vectors $y_{\text{obs}}$, which resulted in
\[
\Sigma_{E_1} = \footnotesize{\left[
\begin{array}{rrrrrrr}
    0.0102  &  0.0006 &  -0.0003 &  -0.0002 &  -0.0001  & -0.0005 &  -0.0000\\
    0.0006  &  0.0093 &  -0.0001  & -0.0004 &   0.0004  & -0.0003 &   0.0004\\
   -0.0003  & -0.0001  &   0.0070  &  0.0002 &  -0.0002  &  0.0002  & -0.0001\\
   -0.0002 &  -0.0004  & 0.0002 &   0.0089 &   0.0002 &  -0.0001 &   0.0002\\
   -0.0001  &  0.0004 &  -0.0002 &   0.0002 &   0.0087 &  -0.0001 &   0.0002\\
   -0.0005  & -0.0003 &   0.0002  & -0.0001 &  -0.0001 &   0.0081 &  -0.0003\\
   -0.0000  &  0.0004 &  -0.0001  &  0.0002 &  0.0002 &  -0.0003  &  0.0085
\end{array} \right]}
\]
For the measurement errors we simply took $\Sigma_{E_2}$ as $0.01$ times the unit matrix. Due to the independence assumption, the covariance matrix of the combined error
$E = E_1 + E_2$ equals $\Sigma_{E} = \Sigma_{E_1} + \Sigma_{E_2}$. The probability distribution of the error is then given by \eqref{error}. This results in the likelihood function \eqref{likelihoodGauss}, which is also depicted in Figure \ref{fig:likelihood}.

The maximum likelihood estimate for the parameters is $a = 3.0380, c = 9.7250$.

\emph{The prior and the posterior distribution.} In order not to prejudice the results, we assumed an improper prior; the parameters $a,c$ to be estimated are far from zero, so no cut-off was needed either and we could take the prior as constant everywhere. In this situation, the posterior density is simply proportional to the likelihood function and can be computed by normalization.

The posterior estimators can be obtained by numerical integration. We used a grid of size $(100\times 100)$ on the relevant domain $[2.9,3.2]\times[8.5,11.0]$. The numerical results are summarized in Table \ref{tab:statsposterior}.

\begin{table}[h]
\centering
\begin{tabular}{|c|c|c|c|c|}
\hline
mean $\mu_a$  & standard deviation $\sigma_a$ & mean $\mu_c$ &  standard deviation $\sigma_c$ & correlation $\rho_{ac}$ \\ \hline
3.0380	&  0.0446	& 9.7315 &	0.3923 &	$-0.9808$ \\\hline
\end{tabular}
\caption{Posterior estimators obtained by numerical integration.}
\label{tab:statsposterior}
\end{table}

\begin{figure}[h]
    \centering
    \includegraphics[width = 4cm]{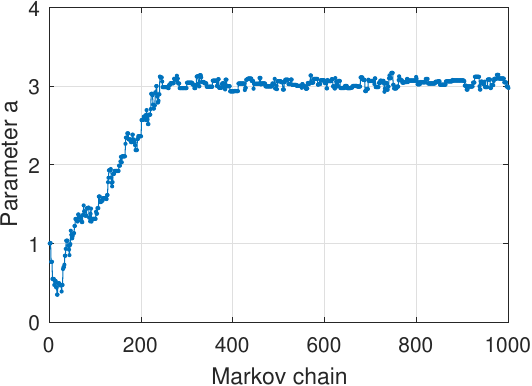} \qquad
    \includegraphics[width = 4cm]{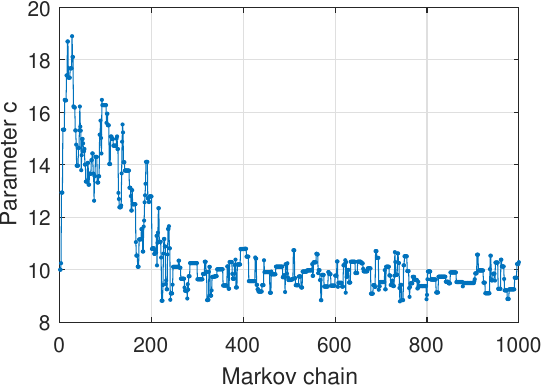}\qquad
    \includegraphics[width = 3.8cm]{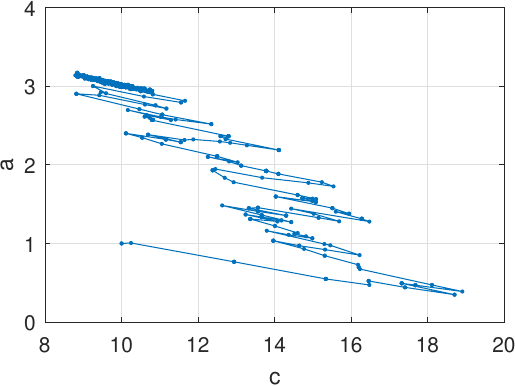}
    \caption{Burn-in phase: first 1000 elements of the Markov chain. Component for $a$ (left), component for $c$ (center), both components (right).}
    \label{fig:burnin}
\end{figure}

\emph{Parameter settings for the Metropolis-Hastings algorithm.} The Markov chain has two components $\xi = (\xi^a,\xi^c)$ corresponding to the parameters $(a,c)$ to be reconstructed. As proposal distribution we chose a two-dimensional Gaussian kernel. We started with a diagonal covariance matrix $\Sigma_{\text{init}}$ with variances $\sigma_a^2 = 0.01, \sigma_c^2 = 1$ corresponding roughly to the mean fluctuations of $(a,c)$ as read off from Figure \ref{fig:likelihood}. For the purpose of demonstration, we took $(\xi^a_0,\xi^c_0) = (1,10)$ as initial point for the chain (one value far from the reference value, one identical with it). We generated chains of length $N = 5000$. A first run of the MCMC-algorithm with $\Sigma_{\text{init}}$ revealed that the components $(\xi^a,\xi^c)$ are highly correlated with a correlation coefficient of $-0.98$. This suggested to implement the correlation in the proposal distribution as well, so that the further MCMC-runs were done with a covariance matrix $\Sigma_{\text{prop}}$ as shown.
\[
   \Sigma_{\text{init}} =  \begin{bmatrix} 0.01 & 0\\ 0 & 1\end{bmatrix}, \qquad
    \Sigma_{\text{prop}} =  \begin{bmatrix} 0.01 & -0.098\\ -0.098 & 1\end{bmatrix}.
\]
The first question to address is the duration of the burn-in phase. Figure \ref{fig:burnin} shows the first 1000 elements of the chain. It can be seen that the burn-in phase ends well before index $k = 300$. We decided to use the chains from $k = 500$ to $k = 5000$ for the subsequent statistical evaluation. The right picture in Figure \ref{fig:burnin} nicely shows the effect of the correlation in the proposal distribution.

\begin{figure}[h]
    \centering
    \includegraphics[width = 4cm]{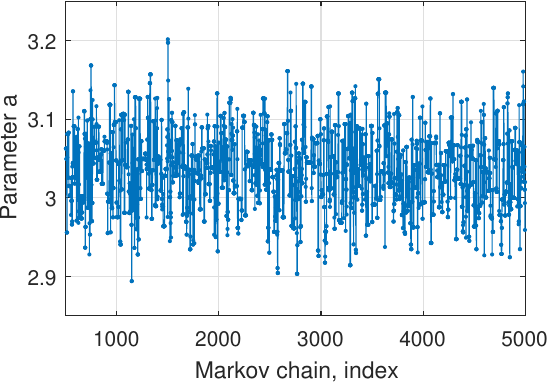} \qquad
    \includegraphics[width = 4cm]{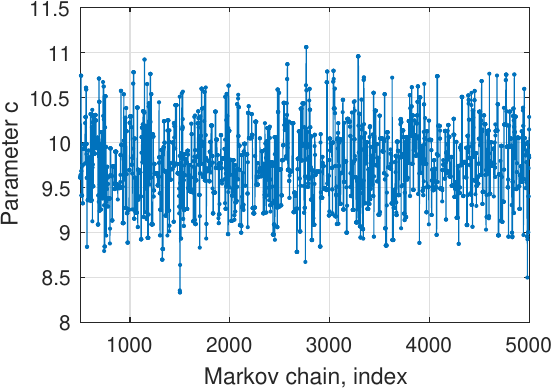}\qquad
    \includegraphics[width = 3.8cm]{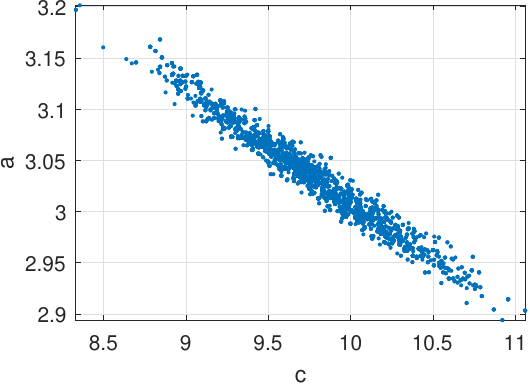}\\[6pt]
    \quad
    \includegraphics[width = 4cm]{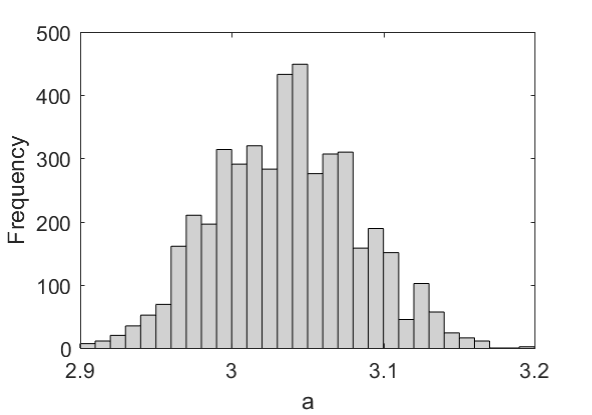} \qquad
    \includegraphics[width = 4cm]{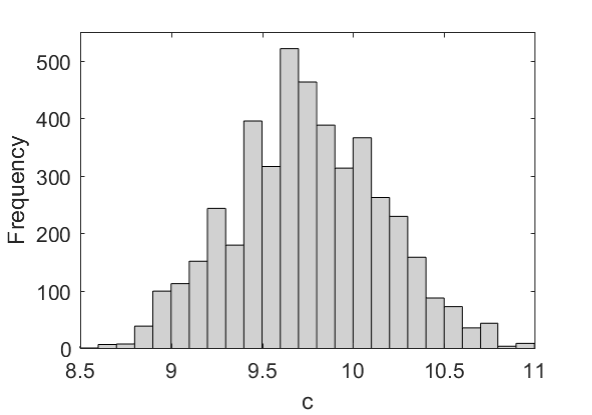}\qquad
    \includegraphics[width = 4.2cm]{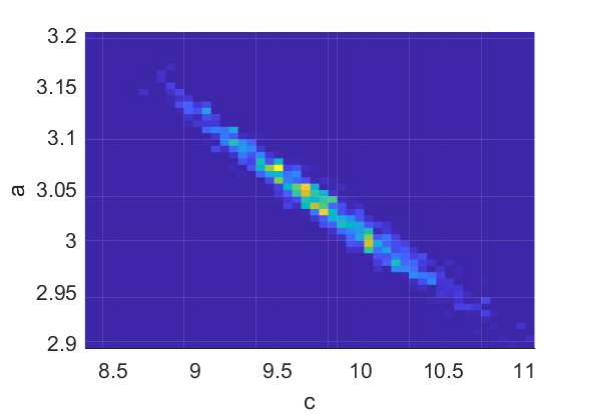}
    \caption{Upper row: Chains after burn-in for parameter $a$ (left), parameter $c$ (center), scatter plot of both parameters (right). Lower row: corresponding histograms. The last figure shows a contour plot of the two-dimensional histogram of $(a,c)$.}
    \label{fig:chains_histograms}
\end{figure}

\emph{Statistical evaluation of the Markov chain.} As noted in Subsection \ref{Metropolis-algorithm}, the end pieces of the Markov chain can serve as a sample of the posterior distribution. The chains from index $k = 500$ to $N = 5000$ are depicted in the upper row of Figure \ref{fig:chains_histograms}. A visual check suggests that the stationary state has been reached and that the mixing properties are good. The acceptance rate turned out to be at the near-optimal value of 24\%.

Histograms of the sample provided by the end pieces of the Markov chain are shown in Figure \ref{fig:chains_histograms}, lower row. The sample statistics may serve as approximations to the corresponding posterior estimators. The results are listed in Table \ref{tab:statsMCMCnoisy}.

\begin{table}[h]
\centering
\begin{tabular}{|c|c|c|c|c|c|c|}
\hline
mean $\mu_a$  & median $m_a$ & std $\sigma_a$ & mean $\mu_c$ & median $m_c$ & std $\sigma_c$ & correlation $\rho_{ac}$\\ \hline
3.0358	& 3.0361 &  0.0468 &	9.7528 & 9.7431 & 0.4085 &	$-0.9816$\\\hline
\end{tabular}
\caption{Posterior estimators obtained from Markov chain, noisy data.}
\label{tab:statsMCMCnoisy}
\end{table}

The somewhat low accuracy is caused by the fact that we used a perturbed trajectory as reference solution. For comparison, we show the same estimates with an unperturbed trajectory as reference solution in Table \ref{tab:statsMCMCexact}.

\begin{table}[H]
\centering
\begin{tabular}{|c|c|c|c|c|c|c|}
\hline
mean $\mu_a$  & median $m_a$ & std $\sigma_a$ & mean $\mu_c$ & median $m_c$ & std $\sigma_c$ & correlation $\rho_{ac}$\\ \hline
3.0024	& 3.0013 &  0.0486 &	9.9900 & 9.9906 & 0.4298 &	$-0.9819$\\\hline
\end{tabular}
\caption{Posterior estimators obtained from Markov chain, exact data.}
\label{tab:statsMCMCexact}
\end{table}

\subsection{An inverse problem for the two-dimensional Poisson equation}
\label{subsec:2dPoisson}

This subsection exemplifies a Bayesian inverse problem in infinite dimensions. Consider the Poisson equation on the square $D = (0,1)\times(0,1)$,
\begin{equation}\label{eq:Poisson}
\begin{array}{rcl}
\Delta u(x_1,x_2) &=& \theta(x_1,x_2),\quad (x_1,x_2) \in D,\\[4pt]
u(x_1,x_2) &=& 0, \quad (x_1,x_2) \in \p D.
\end{array}
\end{equation}
The goal is to reconstruct the right-hand side $\theta$ from an observed solution $u = u_{\text{obs}}$. It is well-known that the forward operator $\cG$ which maps the right-hand side $\theta$ to the solution $u$ is a linear isomorphism from $L^2(D)$ to $H^2(D)\cap H_0^1(D)$. For the purpose of the Bayesian inverse problem, we view $\cG$ as a continuous linear operator from $\cH = L^2(D)$ to $\cH' = L^2(D)$.

As examples, we will take right-hand sides (a) of pyramid shape and (b) equal to a constant, namely
\[
\text{(a)}\quad \theta(x_1,x_2) = 0.1 + 0.9\big(\tfrac12 - |x_1-\tfrac12|\big)\big(\tfrac12 - |x_2-\tfrac12|\big), \qquad\qquad \text{(b)}\quad \theta(x_1,x_2)\equiv 1.
\]
Equation \eqref{eq:Poisson} will be solved numerically by a standard finite difference scheme as described e.g. in \cite[Section 11.2]{Burden:1978} on a grid of size $20\times 20$.

\emph{Observed data.} The observed data are artificially generated by perturbing the right-hand side with independent Gaussian noises at the interior grid points, that is, we add a random realization of a Gaussian $\cN(0,\sigma^2_{\text{noise}}\,I_{19\times 19})$ to the discretized right-hand sides, choosing
\[
\text{(a)}\quad \sigma_{\text{noise}} = 0.01, \qquad\qquad \text{(b)}\quad \sigma_{\text{noise}} = 0.1.
\]
The corresponding (noisy) numerical solutions $u\in L^2(D)$ of \eqref{eq:Poisson} serve as observations. For the sake of consistency with the previous notation, we write again
$y = y_{\rm obs}$ for the observed data.

\emph{The likelihood function.} We take the likelihood function of the form \eqref{eq:likelihoodinfdim}, that is,
\begin{equation*}
\ell(\theta,y) = \exp\Big(-\frac1{2\sigma^2}\|y - \cG(\theta)\|_{L^2(D)}^2\Big)
\end{equation*}
where the error variance $\sigma^2$ is considered as a parameter to be tuned in the algorithm. In the discretization, the $L^2(D)$-norm is approximated by Riemann sums.

\emph{The prior distribution and the solution of the inverse problem.} The priors are taken as centered Gaussian measures $\pi_0(\dd\theta)\sim\cN(0,\cC_0)$ on $L^2(D)$. The covariance operator is taken of the form $\cC_0 = \tau\cC_\alpha$ with $\cC_\alpha$ as defined in Example 3 of Subsection \ref{subsec:Gaussianmeasures} and a scaling parameter
$\tau > 0$.
As indicated in Subsection \ref{subsec:inverseHilbert}, the solution of the Bayesian inverse problem is considered to be the measure $\pi^y(\dd\theta)$ given by \eqref{eq:solutionIP}. The reconstruction of the right-hand side of \eqref{eq:Poisson} is obtained as the mean of a realization of a random field distributed according to
$\pi^y(\dd\theta)$. This realization is obtained using the pCN-algorithm described in Subsection \ref{subsec:MCMCfunctions}.

The starting values $\xi_0$ of the chains are chosen as
\[
\text{(a)}\quad \xi_0(x_1,x_2) = \big(\tfrac12 - |x_1-\tfrac12|\big)\big(\tfrac12 - |x_2-\tfrac12|\big), \qquad\qquad \text{(b)}\quad \xi_0\equiv 0.
\]
Experience from the finite dimensional example in Subsection \ref{subsec:nonlinearspring} shows that the starting value should be chosen not too far away from the parameters to be reconstructed (to avoid a long burn-in phase), whence the choice in (a). In case (b) it turned out that the starting value need not be close to the constant to be reconstructed. To demonstrate this, we chose zero as initial value.

\emph{Parameter settings for the pCN-algorithm.} There are four parameters in the algorithm: the two variance scaling parameters $\sigma,\tau$, the parameter $\beta$ in the proposal step, and the smoothness parameter $\alpha$ of the prior covariance operator. In addition, the Karhunen-Lo\`{e}ve expansion needs to be truncated at length $N_{\rm KL}$. Further, it has to be checked what length $N_{\rm MC }$ of the Markov chain is required and what length $N_{\rm BI}$ of the burn in phase should be chosen. In all cases, we choose
\[
   N_{\rm KL} = 500, \quad N_{\rm MC} = 5000, \quad N_{\rm BI} = 2500.
\]
We tried various values of the smoothness parameter $\alpha$. From the theory, we must take $\alpha > 1$. The larger $\alpha$, the faster is the convergence of the Karhunen- Lo\`eve expansion and the smoother is the corresponding random field. Generally, $\alpha = 4$ turned out to be a good choice, but we also considered $\alpha = 3$, $\alpha = 2$.
For the remaining parameters one can find some recommendations in the literature \cite[Section 4.2]{Cotter:2013}, but in the end one has to tune them experimentally by trial and error.

Our final choices were
\[
\text{(a)}\quad \beta = 0.02,\ \sigma = 0.001,\ \tau = 10, \qquad\qquad \text{(b)}\quad \beta = 0.02,\ \sigma = 0.0001,\ \tau = 1.
\]

\begin{figure}[h]
    \centering
    \includegraphics[width = 5cm]{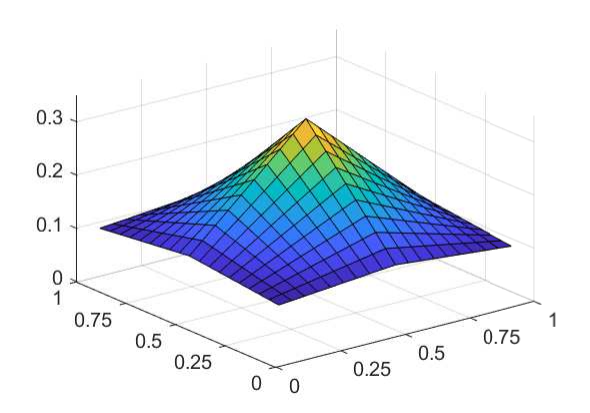} \quad
    \includegraphics[width = 5cm]{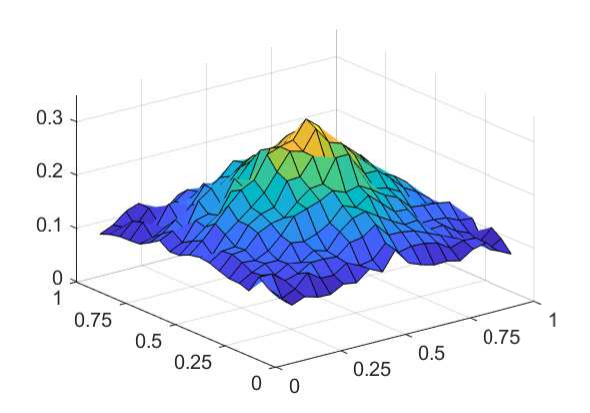}\quad
    \includegraphics[width = 5cm]{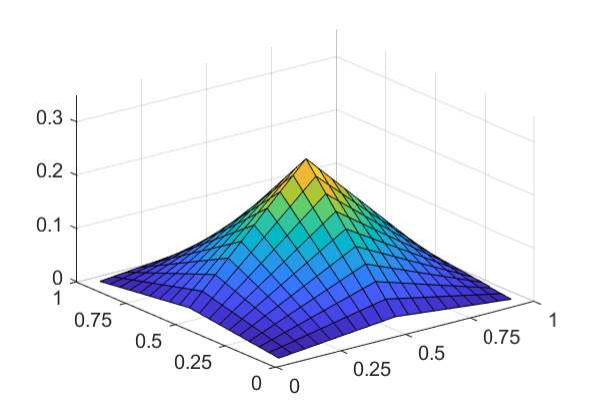}\\[6pt]
    \includegraphics[width = 5cm]{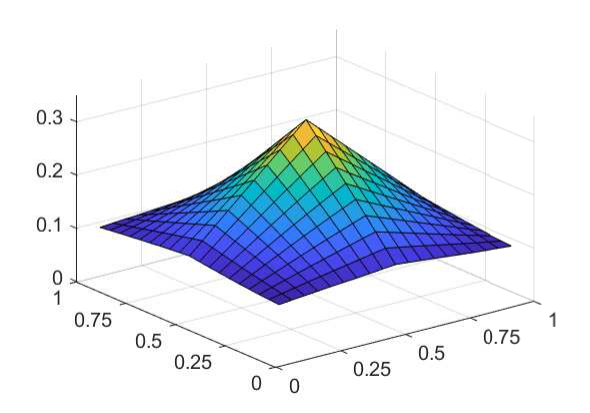} \qquad
    \includegraphics[width = 5cm]{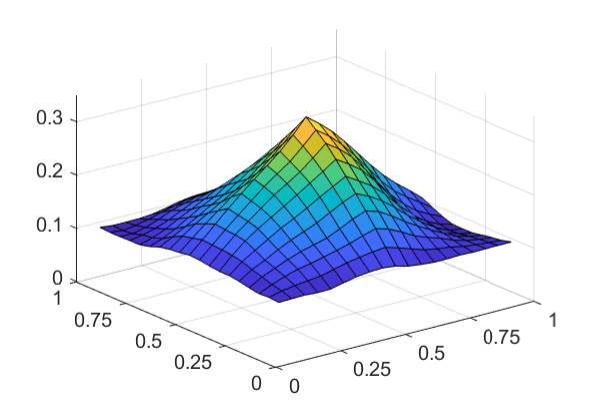}
    \caption{Case (a): reference right-hand side (top left), noisy perturbation of right-hand side (top center), starting value (top right); reconstruction with $\alpha = 4$ (bottom left), reconstruction with $\alpha = 3$ (bottom right). Note that the surface plots show only the interior points of the square, as computed by the finite difference scheme. The missing boundary data are constants.}
    \label{fig:pyramid_surfaces}
\end{figure}

\emph{Results for case (a).} Figure \ref{fig:pyramid_surfaces} shows the reference right-hand side to be reconstructed, the noisy version used for computing the observed data, the chosen starting value and the reconstructed right-hand side with $\alpha = 4$ and $\alpha = 3$. One can see that the reconstruction with $\alpha = 3$ exhibits some random perturbations. The acceptance rates were about 15\% in both cases.

\begin{figure}[h]
    \centering
    \includegraphics[width = 5cm]{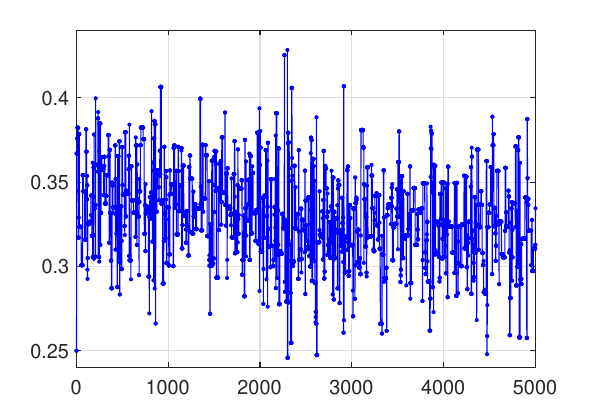} \quad
    \includegraphics[width = 5cm]{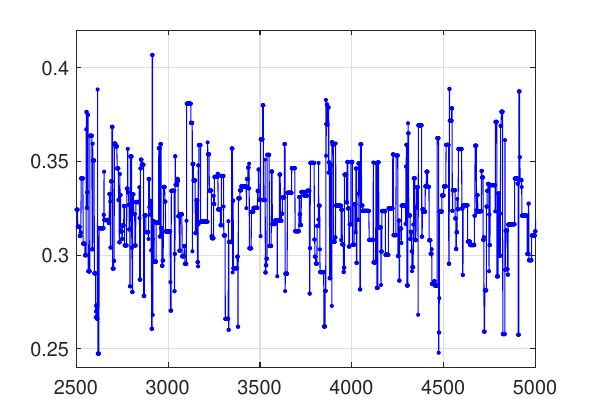}\quad
    \includegraphics[width = 5cm]{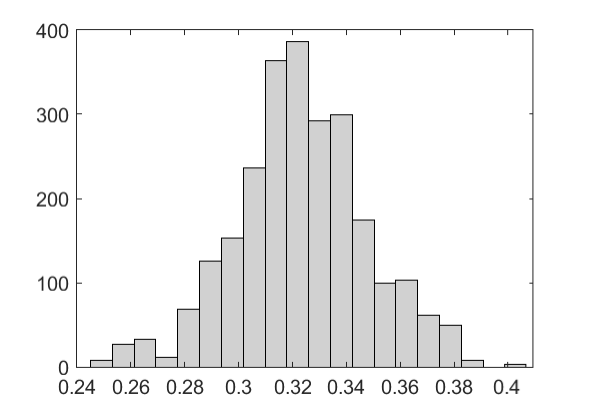}
    \caption{Case (a): Markov chain at center point (left), Markov chain after burn-in (center), histogram of values sampled at center point after burn-in (right).}
    \label{fig:pyramid_chains_histograms}
\end{figure}

The pCN-algorithm produces $19^2$ individual (correlated) chains at each interior grid point. For further illumination, we show the behavior of the chain corresponding to the center point $x_1 = x_2 = 1/2$. Figure \ref {fig:pyramid_chains_histograms} shows the chain including burn-in phase, the second half of the chain, and a histogram of the random realization of $\Theta$ at this point with $\alpha = 4$, produced from the second half of the chain. Table \ref{tab:statscenter} lists some statistical estimates of $\Theta$ evaluated at the center point.

\begin{table}[h]
\centering
\begin{tabular}{|c|c|c|c|c|}
\hline
actual reference value & mean for $\alpha = 4$ & std for $\alpha = 4$ &  mean for $\alpha = 3$ & std for $\alpha = 3$ \\ \hline
0.325	&  0.3234	& 0.0246 &	0.3285 &	0.0270 \\\hline
\end{tabular}
\caption{Case (a): Statistical evaluation of Markov chain at center point after burn-in.}
\label{tab:statscenter}
\end{table}

\begin{figure}[h]
    \centering
    \includegraphics[width = 5cm]{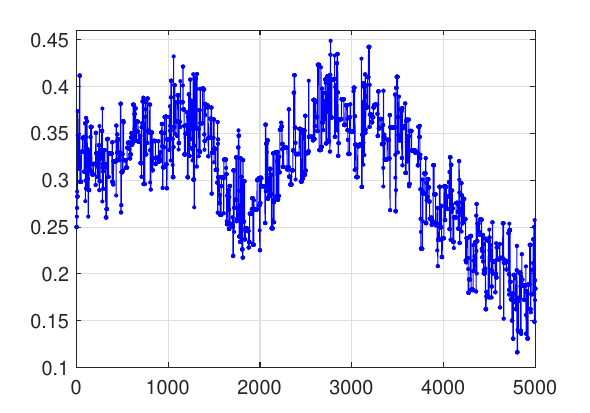} \quad
    \includegraphics[width = 5cm]{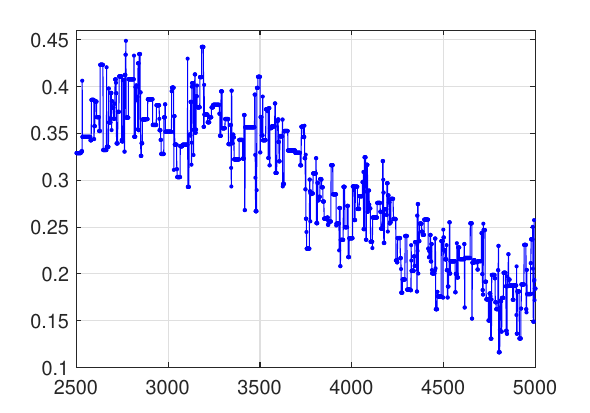}\quad
    \includegraphics[width = 5cm]{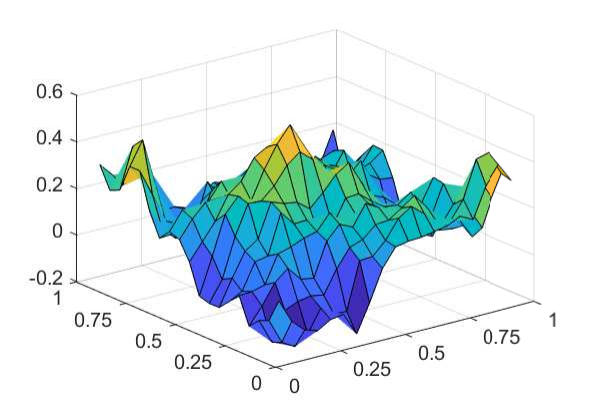}
    \caption{Case (a): Failure of convergence when $\alpha = 2$. Chain at center point (left), second half of chain at center point (center), failed reconstruction of right-hand side (right).}
    \label{fig:pyramid_alpha_2}
\end{figure}

For completeness, we show that for $\alpha = 2$ the chain fails to converge (at least up to $N_{\rm MC} = 5000$). Figure \ref{fig:pyramid_alpha_2} depicts the chain evaluated at the center point and the reconstructed mean. The bad behavior of the chain is clearly visible; its mean after burn-in equals $0.2941$, its standard deviation is $0.0748$.
The reconstruction of the right-hand side is equally wrong.

\begin{figure}[H]
    \centering
    \includegraphics[width = 5cm]{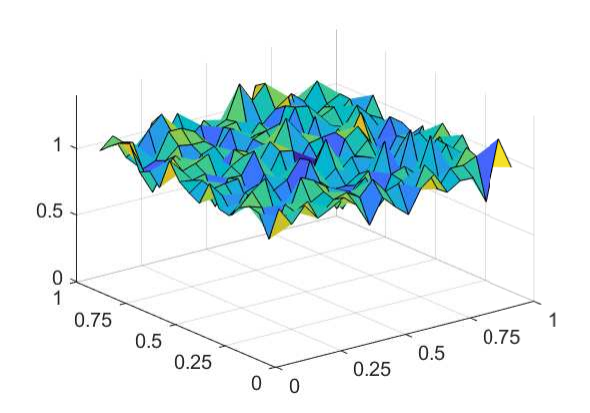} \quad
    \includegraphics[width = 5cm]{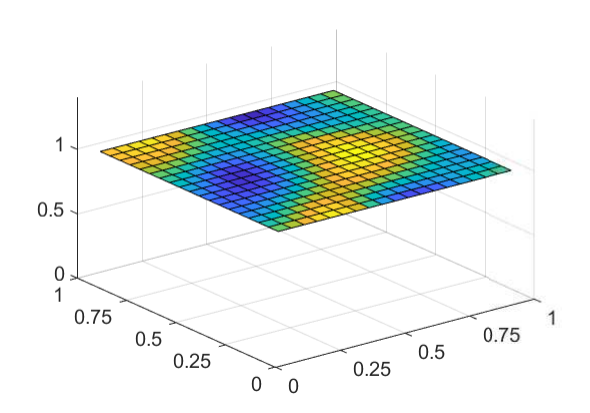}\quad
    \includegraphics[width = 5cm]{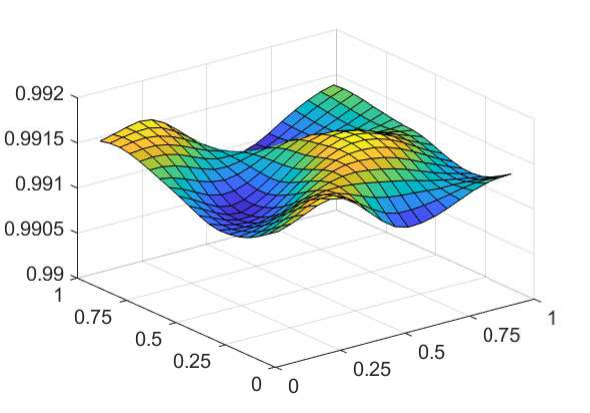}
    \caption{Case (b): Constant reference right-hand side perturbed by noise (left), reconstructed constant (center), zoom near height one to visualize the small fluctuation of the reconstructed constant (right).}
    \label{fig:constant_surfaces}
\end{figure}

\begin{figure}[H]
    \centering
    \includegraphics[width = 5cm]{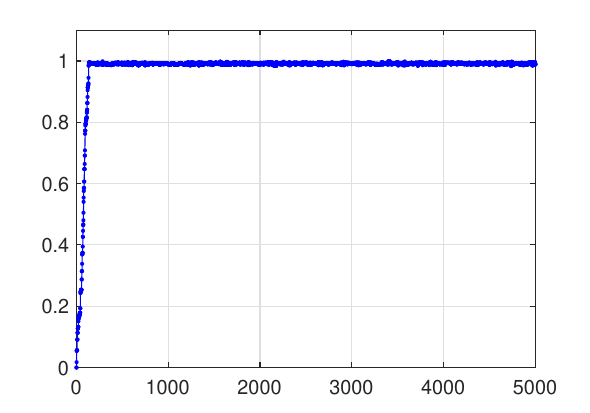} \quad
    \includegraphics[width = 5cm]{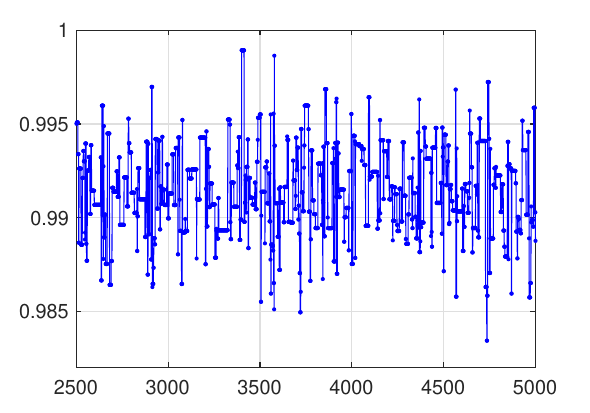}\quad
    \includegraphics[width = 5cm]{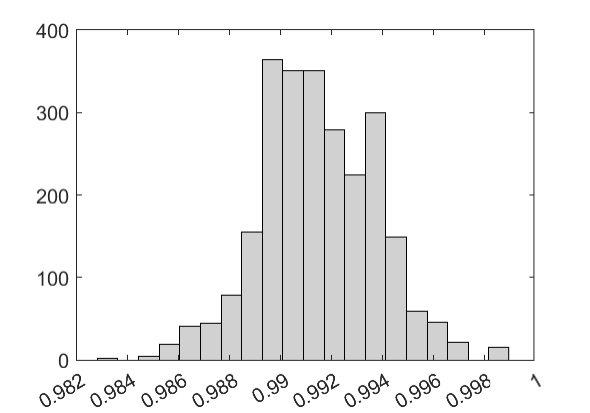}
    \caption{Case (b): Markov chain at center point (left), Markov chain after burn-in (center), histogram of values sampled at center point after burn-in (right).}
    \label{fig:constant_chains_histograms}
\end{figure}

\emph{Results for case (b).} In this case, we seek to reconstruct a constant right-hand side (equal to one). As initial value of the chain we took $\xi_0 \equiv 0$. The acceptance rate was 16.5\%. Figure \ref{fig:constant_surfaces} shows the noisy reference right-hand side and the reconstructed constant, choosing $\alpha = 4$. The reconstruction is very smooth; in order to make the small fluctuations visible, we also zoom in near height one. In Figure \ref{fig:constant_chains_histograms}, the behavior of the chain at the center point $x_1 = x_2 = 1/2$ as well as a histogram of the random sample of $\Theta$ at this point (generated from the second half of the chain) can be seen.
Since the reconstructed right-hand side is (nearly) constant, it makes sense to evaluate the total mean and standard deviation of all $19^2$ chains at the interior grid points.
The values are listed in Table \ref{tab:statscenterconst}, together with the values obtained from observations generated by the unperturbed, strictly constant right-hand side.
We also executed the pCN-algorithm with $\alpha < 4$. We found convergence down to $\alpha = 2$, but not below.
\begin{table}[H]
\centering
\begin{tabular}{|c|c|c|c|}
\hline
total mean, noisy data & total std, noisy data &  total mean, exact data & total std, exact data \\ \hline
 0.9913	& 0.000159 &	0.999970 &	0.000135 \\ \hline
\end{tabular}
\caption{Case (b): Statistical evaluation of Markov chain at center point after burn-in, noisy data vs. exact data.}
\label{tab:statscenterconst}
\end{table}

\section{Tools from probability and measure theory}
\label{Sec:appendix}
\subsection{Conditional distributions}
\label{subsec:probability}
In the finite dimensional case treated in Section \ref{sec:Bayes}, the elementary notion of conditional probability densities $\pi(\theta|y)$ turned out to be sufficient. However, for the theory of Markov chains and for Bayesian inversion in infinite dimensional spaces, the general notion of conditional distributions of random variables is required. The goal of this subsection is to rigorously introduce the concepts of conditional expectation, conditional probability and conditional distributions. A short introduction to conditional probabilities is in \cite[Section 1.7]{Arnold:1974}; full details can be found in \cite[Chapter X]{Bauer:1968} and in \cite[Sections 10.1 and 10.2]{Dudley:2002}.

In what follows ${\mathcal H}, {\mathcal H}'$ will be separable Hilbert spaces, equipped with their Borel $\sigma$-algebras ${\mathcal B} = {\mathcal B}({\mathcal H})$,
${\mathcal B}' = {\mathcal B}({\mathcal H}')$. Let $(\Omega,{\mathcal A},P)$ be a probability space, ${\mathcal A}'$ a sub-$\sigma$-algebra of ${\mathcal A}$ and let $\Theta \in L^1(\Omega)$ be an $\cA$-measurable and $P$-integrable random variable with values in ${\mathcal H}$. In general, $\Theta$ is no longer measurable with respect to the coarser
$\sigma$-algebra $\cA'$. However, there is an (almost surely unique) $\cA'$-measurable random variable $\Theta'$ such that
\[
\int_{A'}\Theta'(\omega)P({\rm d}\omega) = \int_{A'}\Theta(\omega)P({\rm d}\omega)
\]
for all $A'\in{\mathcal A}'$. Its existence follows from an application of the Radon-Nikodym theorem. The random variable $\Theta'$ is called the  \emph{conditional expectation of $\Theta$ given ${\mathcal A}'$} and is denoted by ${\rm E}(\Theta|{\mathcal A}')$; we have that
\[
\int_{A'}{\rm E}(\Theta|{\mathcal A}')(\omega)P({\rm d}\omega) = \int_{A'}\Theta(\omega)P({\rm d}\omega)
\]
for all $A'\in{\mathcal A}'$. If $\Theta\in L^2(\Omega)$, the conditional expectation can also be viewed as the orthogonal projection of $\Theta$ to the subspace of $L^2(\Omega)$ of $\cA'$-measurable functions \cite[Theorem 15.8]{Bauer:1996}.

Given $A\in\cA$, we can apply the construction of conditional expectation to the real-valued random variable $\eins_A$. The \emph{conditional probability of the event $A\in {\mathcal A}$ given ${\mathcal A}'$} is defined as
$P(A|{\mathcal A}')(\omega) = E(\eins_A|{\mathcal A}')(\omega)$. It has the property that
\[
  \int_{A'}P(A|{\mathcal A}')(\omega)P({\rm d}\omega) = P(A'\cap A), \quad A\in{\mathcal A}, A'\in{\mathcal A}'.
\]
\emph{Example.} If ${\mathcal A}' = \{\emptyset, A', A'' =\Omega\setminus A',\Omega\}$, then $P(A|{\mathcal A}')(\omega) = P(A'\cap A)/P(A')$ for $\omega\in A'$ and similarly with $A''$ in the role of $A'$.

Let $Y: \Omega\to\cH'$ be an $\cA$-$\cB'$ measurable function and let $\cA'$ be the $\sigma$-algebra generated by the inverse images $Y\inv(B')$, $B'\in\cB'$. The \emph{conditional expectation of $\Theta$ given $Y$} is $\rE(\Theta|Y) = \rE(\Theta|\cA')$. Similarly, the conditional probability of an event $A\in\cA$ given $Y$ is $P(A|Y) = \rE(\eins_A|\cA')$.
For $A\in\cA$, the conditional probability $P(A|Y)$ is a real valued, $\cA'$-measurable function. By the factorization property \cite[Lemma 55.1]{Bauer:1968}, there is a real-valued function $h:\cH'\to\R$, measurable with respect to $\cA'$, such that
\[
   P(A|Y)(\omega) = h(Y(\omega)).
\]
This gives rise to the definition $P(A|Y=y) = h(y)$. Note that $P(A|Y)(\omega)$ is constant on $\{\omega: Y(\omega) = y\}$, and we have that
$P(A|Y)(\omega) = P(A|Y=Y(\omega))$.

Applying this construction to $A = \Theta\inv(B)$, $B\in\cB$ leads to the notion of \emph{conditional distribution of $\Theta$ given $Y$}. This is understood as a probability kernel $p(\omega,B)$ such that
\begin{itemize}
\item[(a)] for all $\omega\in\Omega$, $p(\omega,\cdot)$ is a probability measure on $\cH$;
\item[(b)] for all $B\in\cB$, $p(\cdot,B)$ is a version of $P(\Theta\inv(B)|Y)$.
\end{itemize}
For $B\in\cB$, this means that $p(\omega,B) = P(\Theta\inv(B)|Y)(\omega)$ for almost all $\omega$. The exceptional set of measure zero a priori depends on $B$. The fact that a probability kernel exists is nontrivial and uses the property that $\cH$ is a Polish space (i. e. a topological space with a countable base and topology generated by a complete metric), see e. g. \cite[Satz 56.5]{Bauer:1968}, \cite[Theorem 10.2.2]{Dudley:2002}.

As above, we can interpret $P(\Theta\inv(B)|Y)(\omega)$ as a function of $y\in\cH'$, which we denote by
\[
   \pi^y(B) = P(\Theta\inv(B)|Y=y).
\]
Note that $p(\omega,B) = \pi^{Y(\omega)}(B)$. At fixed $y = Y(\omega)$, $B\to \pi^y(B)$ is a probability measure on $\cH$. We denote this probability measure by $\pi^y(\dd\theta)$ and also refer to it as the conditional distribution of $\Theta$ given $y$. (The notation follows \cite{Stuart_2010_acta_numerica} and is widely used in the literature on Bayesian inverse problems.)

In the finite dimensional case, where $\cH=\R^d$, $\cH'=\R^{d'}$ and the random variable $(\Theta,Y)$ has a probability density $\pi(\theta,y)$ with respect to Lebesgue measure $\dd\theta\times \dd y$, it can be shown \cite[Theorem 15.11]{Bauer:1996} that $\pi^y(\dd\theta)$ has a probability density as well, and it coincides with the conditional density $\pi(\theta|y)$.

\subsection{Markov chains}
\label{subsec:Markovchains}

In order to analyze the behavior of MCMC-algorithms, basic notions from the theory of Markov chains are needed. In this part of Section \ref{Sec:appendix}, we collect a minimal selection of required concepts; we base the exposition on \cite{robert1999monte}. We also refer to the readable survey \cite{Tierney:1994} and to \cite{MeynTweedie:1993,Nummelin:1984}.

A Markov chain is a sequence of random variables $X_0, X_1, X_2,\ldots$ on a probability space $(\Omega,\cA,P)$ with values in a measurable space $(S,\cS)$, the \emph{state space}. The Markov property says that, for all $n\in\N$, the conditional probability distributions satisfy
\[
P(X_{n+1}|X_0, X_1,\ldots, X_n) = P(X_{n+1}|X_n).
\]
We assume that the Markov chain is time homogeneous, that is $P(X_{n+1}|X_n) = P(X_1|X_0)$ for all $n\in\N$. A \emph{transition kernel} $p(x,A)$ is a function on $S\times\cS$ such that
\begin{itemize}
\item[(a)] for all $x\in S$, $p(x,\cdot)$ is a probability measure;
\item[(b)] for all $A\in\cS$, $p(\cdot,A)$ is measurable.
\end{itemize}
The transitions of the chain are determined by the kernel:
\[
  P(X_{n+1}\in A|X_n= x) = p(x,A) = \int_A p(x,\dd y)
\]
for all $x\in S$, $A\in \cS$. As noted in Subsection \ref{subsec:probability}, such a transition kernel exists if $S$ is a Polish space and $\cS = \cB(S)$. (With the notation of
Subsection \ref{subsec:probability}, taking $X_{n+1} = \Theta$, $X_n = Y$, we have $p(x,A) = \pi^x(A)$.)

The finite dimensional joint distributions of the members of the chain are determined by the transition kernel and by specifying the initial distribution, the distribution of $X_0$. If the initial distribution is a Dirac measure $\delta_x$ ($X_0$ has the deterministic value $x$), and $Y$ is a random variable, we write $P(Y\in B|X_0 = x) = P_x(Y\in B)$. Similarly, we write $\rE_x(Y) = \rE(Y|X_0=x)$. We get
\begin{eqnarray}\label{eq:Px}
P_x(X_1\in A_1) &=& p(x,A_1), \nonumber\\
P_x(X_1\in A_1, X_2 \in A_2) &=& \int_{A_1}p(y_1,A_2)\,p(x,\dd y_1),\nonumber\\
   &\vdots &\\
P_x(X_1\in A_1,\ldots, X_n\in A_n) & =& \int_{A_1}\cdots \int_{A_n-1} p(y_{n-1},A_n)\,p(y_{n-2},\dd y_{n-1})\cdots p (x,\dd y_1). \nonumber
\end{eqnarray}
This can be inferred from properties of conditional expectations. see e.g. \cite[Section 2.2]{Arnold:1974}.
Introduce the space of trajectories of the chain as $\Omega^\infty = \prod_{i=0}^\infty \Omega_i$ where each $\Omega_i$ is a copy of $(\Omega,\cA)$ and $\Omega$ is equipped with the product $\sigma$-algebra $\cA^\infty$.
Defining $P_x(A_1\times\cdots\times A_n)$ by the right-hand side of \eqref{eq:Px}, it immediately follows from the defining formula that
$P_x$ can be extended to a measure on $(\Omega^\infty,\cA^\infty)$, see also \cite[Theorem 3.4.1]{MeynTweedie:1993}.

The marginal probability distributions of the chain, conditioned on $X_0=x$, are denoted by
\[
  p^{(n)}(x,A) =  P_x(X_{n}\in A).
\]
(In Subsection \ref{Metropolis-algorithm}, the marginal probability measures were assumed to have densities $p^{(n)}(y)$, suppressing the dependence on the initial value $x$.)

Next, we highlight some important concepts concerning Markov chains.
\begin{itemize}
\item  Let $A\in \cS$. The \emph{number of passages} of $(X_n) = (X_n)_{n\in\N}$ in $A$ is
\[
  \eta_A = \sum_{n=1}^\infty \eins_A(X_n).
\]
\end{itemize}
An important quantity will be the \emph{average number of passages} $\rE_x(\eta_A)$.
\begin{itemize}
\item Let $\varphi$ be a measure on $(S,\cS)$. The chain $(X_n)$ is \emph{$\varphi$-irreducible} if for every $x\in S$ and $A\in\cS$ with $\varphi(A) >0$ it holds that $\rE_x(\eta_A) > 0$ (cf. \cite[Definition 6.13 and Theorem 6.15]{robert1999monte}).
\end{itemize}
\begin{itemize}
\item A $\varphi$-irreducible chain with transition kernel $p$ is periodic if there exists an integer $d\geq 2$ and a sequence $\{S_0,\ldots, S_{d-1}\}$ of $d$ nonempty disjoint sets in $S$ such that, for all $i = 0,\ldots, d-1$ and all $x\in S_i$,
\[
    p(x,S_j) = 1\quad \text{for}\ j = i+1\ (\text{mod}\ d).
\]
Otherwise, the chain is \emph{aperiodic}, cf. \cite[Section 3]{Tierney:1994}.
\end{itemize}
\begin{itemize}
\item A Markov chain $(X_n)$ is \emph{recurrent} if (i) $(X_n)$ is $\varphi$-irreducible for some measure $\varphi$ and (ii) for every $A\in\cS$ such that $\varphi(A) > 0$, $\rE_x(\eta_A) = \infty$ for every $x\in S$.
\end{itemize}
Thus for whatever state the chain is initially, the expected number of visits to the set $A$ is infinite.
A stronger requirement is that the probability of infinite return to $A$ is equal to one:
\begin{itemize}
\item A Markov chain $(X_n)$ is \emph{Harris recurrent} if  $(X_n)$ is $\varphi$-irreducible for some measure $\varphi$ and for all $A\in\cS$ such that $\varphi(A) > 0$,
$P_x(\eta_A = \infty) = 1$ for every $x\in A$.
\end{itemize}
Invariant measures and stationary distribution are central notions of the theory.
\begin{itemize}
\item A $\sigma$-finite measure $\pi$ is invariant for the chain (and its transition kernel) if
\begin{equation}\label{eq:invariantmeasure}
   \pi(B) = \int p(x,B)\pi(\dd x)
\end{equation}
for all $B\in\cS$. If there exists an \emph{invariant probability measure} for a $\varphi$-irreducible chain, the chain is called \emph{positive}. The invariant probability measure is also referred to as \emph{stationary}, because $X_0\sim\pi$ implies $X_n\sim \pi$ for all $n\in\N$; the chain is stationary in distribution. Note that if $\pi$ is a probability measure, the invariance condition \eqref{eq:invariantmeasure} is equivalent with
\begin{equation}\label{eq:invariantprobability}
   \int f(y)\pi(\dd y) = \iint f(y)p(x,\dd y)\pi(\dd x)
\end{equation}
for all bounded measurable functions $f$ on $S$.
\end{itemize}
The question of existence of an invariant measure is answered by the following criterion.
\begin{itemize}
\item Suppose the chain $(X_n)$ is recurrent. Then it has an invariant $\sigma$-finite measure. The measure is unique up to a multiplicative factor.
\end{itemize}
We refer to \cite[Theorem 6.38]{robert1999monte} and to \cite[Section 10]{MeynTweedie:1993} for a proof.
For applications to the Metropolis-Hastings algorithm and the pCN-algorithm, the following simple criterion is useful.
A Markov chain $(X_n)$ with transition kernel $p$ is said to satisfy the \emph{detailed balance condition} if there is a probability measure $\pi$ such that
\begin{equation}\label{eq:detailedbalance}
   p(x,\dd y)\pi(\dd x) = p(y,\dd x)\pi(\dd y)
\end{equation}
as measures on $S\times S$ for all $x, y$.
\begin{proposition}\label{prop:detbal}
If the transition kernel $p$ of a Markov chain $(X_n)$ satisfies the detailed balance condition with a probability measure $\pi$, then $\pi$ is an invariant probability measure of the chain.
\end{proposition}
\emph{Proof.} The assertion is easily proved by the observation that
\[
  \iint f(y)p(x,\dd y)\pi(\dd x) = \iint f(y)p(y,\dd x)\pi(\dd y) = \int f(y)\pi(\dd y)
\]
for every bounded measurable function $f$ on $S$, since $\int \pi(y,\dd x) = 1$ for all y. Thus \eqref{eq:invariantprobability} holds. \qed

From the previous assertion, $\pi$ is unique provided the chain is recurrent.

We are ready to state a first convergence result for Markov chains. The \emph{total variation distance} between two probability measures $\pi_1,\pi_2$ is defined as
\[
\big\|\pi_1 - \pi_2\big\| = \sup_{A\in\cS}|\pi_1(A) - \pi_2(A)|.
\]
We have the following assertion.
\begin{proposition}\label{prop:TV}

Suppose a Markov chain $(X_n)$ with transition kernel $p$ is $\pi$-irreducible with invariant probability measure $\pi$, Harris recurrent and aperiodic. Then
\[
  \lim_{n\to\infty}\left\| \int p^{(n)}(x,\cdot)\mu(\dd x) - \pi\right\| = 0
\]
for every initial distribution $\mu$. In particular,
\[
  \lim_{n\to\infty}\left\| p^{(n)}(x,\cdot) - \pi\right\| = 0
\]
for every $x\in S$.
\end{proposition}
\emph{Proof.} The result is stated in \cite[Theorem 6.51]{robert1999monte}, \cite[Theorem 1]{Tierney:1994}. For a proof, see \cite[Corollary 6.7(ii)]{Nummelin:1984}. \qed

A chain which satisfies the hypotheses in Proposition \ref{prop:TV} is often called \emph{ergodic}. One can obtain convergence rates under stronger conditions, such as geometric or uniform ergodicity \cite[Sections 6.6.2 and 6.6.3]{robert1999monte}.
Note that Proposition \ref{prop:TV} implies that
\[
   \lim_{n\to\infty}p^{(n)}(x,A) = \pi(A)
\]
for every $A\in\cS$ and $x\in S$. Therefore, $\pi$ is often called the \emph{equilibrium distribution} of the chain, and the chain is said to approach its \emph{stationary state} as $n\to\infty$.

The next important issue is the convergence of averages in the sense of a law of large numbers. In the context of Markov chains, this result is usually referred to as the \emph{ergodic theorem}.
\begin{proposition}\label{prop:ET}
Let the Markov chain $(X_n)$ be $\pi$-irreducible with invariant probability measure $\pi$ and Harris recurrent. Let $f(x)$ be any function integrable with respect to $\pi(\dd x)$. Then
\[
  \lim_{N\to\infty} \frac1{N} \sum_{n=1}^N f(X_n) = \int_S f(x)\pi(\dd x),
\]
$P_{x_0}$-almost surely for every initial value $x_0$, where $P_{x_0}$ is the measure on the space of trajectories introduced after \eqref{eq:Px}.
\end{proposition}
\emph{Proof.} We refer to \cite[Theorem 17.1.7]{MeynTweedie:1993} and \cite[Theorem 6.63]{robert1999monte}.\qed

The ergodic theorem is taken as a justification of Monte Carlo sampling using the chain. In fact, let $X$ be a random variable distributed according to the equilibrium distribution $\pi$ and let $\xi_0,\xi_1,\xi_2,\ldots $ be a realization of the chain.
Then Proposition \ref{prop:ET} suggests that
\[
   \lim_{N\to\infty} \frac1{N} \sum_{n=1}^N f(\xi_n) = \rE(f(X)).
\]
For sufficiently large $M$ (when the chain is close to its stationary state), the end pieces $\xi_M, \dots, \xi_N$ can be treated as a sample of $X$. In particular,
the sample mean, moments and quantiles of $\xi_M, \dots, \xi_N$ may serve as estimators for the expectation value, moments and quantiles of $X$.

The final issue addresses \emph{central limit theorems for Markov chains.} There are many versions which usually require so-called mixing conditions which are difficult to verify, see the discussion in \cite{Robert:1994}. We focus on one of the versions. Let $\cA,\cA'$ be sub-$\sigma$-algebras of $\cS$. We recall the following measures of dependence:
\begin{itemize}
\item[] $\alpha(\cA,\cA') = \sup|P(A\cup A') - P(A)P(A')|,\ A \in\cA, A'\in \cA'$
\item[] $\varphi(\cA,\cA') = \sup|P(A'|A) - P(A')|,\ A \in\cA, A'\in \cA', P(A) > 0$.
\end{itemize}
We have that $2\alpha(\cA,\cA') \leq \varphi(\cA,\cA')$ \cite[Formula (1.11)]{Bradley:2005}.
Let $(X_n)$ be a stationary (i.e., shift invariant in distribution) Markov chain. Define $\cA_j$ as the sub-$\sigma$-algebra of $\cS$ generated by $X_j$. Put
\[
   \alpha(n) = \alpha(\cA_0,\cA_n),\quad \varphi(n) = \varphi(\cA_0,\cA_n).
\]
(This is a specialization of the general definition to Markov chains \cite[Section 3.1]{Bradley:2005}.) The chain $(X_n)$ is called \emph{$\alpha$-mixing}, if $\alpha(n)\to 0$ as $n\to\infty$, and similarly for \emph{$\varphi$-mixing}. Here is one of the limit theorems.
\begin{proposition}\label{prop:CL}
Suppose that the chain $(X_n)$ is stationary and $\alpha$-mixing with $\alpha(n) = \cO(n^{-5})$ and that $\rE(X_n) = 0$ and $\rE(X_n^{12}) < \infty$. Let
$S_N = X_1 + \cdots + X_N$. Then
\begin{equation}\label{eq:gamma2}
   \lim_{N\to\infty}\frac{\rV(S_N)}{N} = \gamma^2 = {\rm V}(X_1) + 2\sum_{k=1}^\infty{\rm COV}(X_1,X_{1+k})
\end{equation}
where the series is absolutely convergent. If $\gamma > 0$, the random variable $S_N/\gamma\sqrt{N}$ converges (in distribution) to the standard normal distribution.
\end{proposition}
\emph{Proof.} The theorem is taken from \cite[Theorem 27.4]{Billingsley:1995}. \qed

It is stated in \cite{Billingsley:1995} that the conditions are stronger than necessary. Also, there are many results on the decay of mixing coefficients. One of the results is
\cite[Theorem 4 and discussion before]{Davydov:1974}: Assume that $(X_n)$ has an invariant probability distribution, is Harris recurrent and aperiodic. If for some $n_0$, $\varphi(n_0) < 1/2$, then $\varphi(n) \leq C\rho^n$
for some $C>0$, $0 < \rho < 1$ and all $n$. Since $2\alpha(n)\leq \varphi(n)$, the same then holds for $\alpha(n)$, so that the first condition in Proposition \ref{prop:CL} is satisfied; see also \cite[Theorem 3.4]{Bradley:2005}.

\subsection{Convergence of MCMC-algorithms}
\label{subsec:MCMCconvergence}

In this Subsection, we collect results from the literature that imply the convergence of the Metropolis-Hastings algorithm (Subsection \ref{Metropolis-algorithm}) and the pCN-algorithm (Subsection \ref{subsec:MCMCfunctions}).
\smallskip

\emph{The Metropolis-Hastings algorithm.} The chain produced by the algorithm is a Markov chain which has a transition kernel $p(x,A)$.
Sampling from the proposal distribution $q(x,y)$ together with the accept/reject step with acceptance probability $a(x,y)$ determine the kernel as
\[
   p(x,A) = \int a(x,y)q(x,y)\dd y + (1-b(x)) \delta_x(A) \quad\text{with}\quad b(x) = \int a(x,y)q(x,y)\dd y
\]
where $\delta_x$ is the Dirac mass in $x$, or
\begin{equation}\label{eq:MCMCkernel}
  p(x,\dd y) = a(x,y)q(x,y)\dd y + (1-b(x))\delta_x(\dd y).
\end{equation}
We are going to show that the kernel $p$ satisfies the detailed balance condition \eqref{eq:detailedbalance} with the probability measure $\pi$. Recalling that $\pi(\dd x) = \pi(x)\dd x$ and integrating the first part against a bounded continuous function $f(x,y)$ we get
\[
   \iint f(x,y)a(x,y)q(x,y)\pi(x)\dd y\dd x =  \iint f(x,y)a(y,x)q(y,x)\pi(y)\dd x\dd y
\]
where we have used the symmetry of $q$ and the fact that $\min(1,\frac{\beta}{\alpha})\alpha = \min(1,\frac{\alpha}{\beta})\beta$. Concerning the second part, we have
\begin{eqnarray} \label{eq:balancedelta}
\iint f(x,y) (1-b(x))\delta_x(\dd y)\pi(x)\dd x &=& \int f(x,x) \pi(x)\dd x = \int f(y,y) \pi(y)\dd y\\
      &=& \iint f(x,y) (1-b(y))\delta_y(\dd x)\pi(y)\dd y. \nonumber
\end{eqnarray}
Combining the two calculations, we get
\[
   \iint f(x,y)p(x,\dd y)\pi(x)\dd x = \iint f(x,y)p(y,\dd x)\pi(y)\dd y
\]
thus the detailed balance condition holds, and $\pi(\dd x) = \pi(x)\dd x$ is an invariant probability measure for the chain.

To prove convergence, we rely on results from \cite{robert1999monte}. First, \cite[Lemma 7.3]{robert1999monte} says that:
\begin{itemize}
\item If the Metropolis-Hastings chain is $\pi$-irreducible, it is Harris recurrent.
\end{itemize}
In particular, $\pi$ is the unique invariant probability measure of the chain. A central criterion for convergence is as follows \cite[Lemma 7.6]{robert1999monte},
\cite[Theorem 2.2]{RobertsTweedie:1996}:
\begin{itemize}
\item Assume that the probability density $\pi(x)$ is bounded away from 0 and $\infty$ on every compact subset of its support. If there exist positive numbers $\eps$ and $\delta$ such that
\[
  q(x,y) > \eps \quad\text{if}\quad |x-y| < \delta,
\]
then the Metropolis-Hastings chain is $\pi$-irreducible and aperiodic.
\end{itemize}
In conclusion, under the conditions above, the Metropolis-Hastings chain has $\pi$ as its invariant probability distribution, is $\pi$-irreducible, Harris recurrent and aperiodic. Thus it satisfies all requirements leading to convergence and ergodicity (Propositions \ref{prop:TV} and \ref{prop:ET}).
\smallskip

\emph{The pCN-algorithm.} To prove that the transition kernel $p$ corresponding to the pCN-algorithm satisfies the detailed balance condition \eqref{eq:detailedbalance} with measure $\pi$, we confine ourselves to quoting the arguments from \cite{DashtiStuart:2017Bayesian}. For the sake of consistency, we denote the elements of the state space $\cH$ by $x, y$ here.
Similarly to \eqref{eq:MCMCkernel}, the transition kernel is
\begin{equation*}
  p(x,\dd y) = a(x,y)q(x,\dd y) + (1-b(x))\delta_x(\dd y) \quad\text{with}\quad b(x) = \int a(x,y)q(x,\dd y)
\end{equation*}
where $q(x,\dd y)$ is the kernel of the proposal distribution, and we have to show that
\begin{equation*}
   p(x,\dd y)\pi(\dd x) = p(y,\dd x)\pi(\dd y).
\end{equation*}
As seen by the argument in \eqref{eq:balancedelta}, the term with the Dirac mass is symmetric in $(x,y)$, thus it suffices to show that
\begin{equation}\label{eq:detailedbalancepCNq}
   a(x,y)q(x,\dd y)\pi(\dd x) = a(y,x)q(y,\dd x)\pi(\dd y).
\end{equation}
To prove this, one first shows that $q$ itself satisfies the detailed balance condition with respect to $\pi_0$,
\begin{equation}\label{eq:auxbalance}
q(x,\dd y)\pi_0(\dd x) = q(y,\dd x)\pi_0(\dd y).
\end{equation}
The form of the kernel $q(x,\dd y)$ is suggested by the proposal step in the pCN-algorithm. Namely, for fixed $x$, $q(x,\dd y)$ is given by the measure $\cN(\sqrt{1-\beta^2}\,x, \beta^2\cC_0)$. Note that a random variable $v$ is distributed according to $q(x,\dd y)$ if $v = \sqrt{1-\beta^2}\,x + \beta w$ where
$w\sim \cN(0,\cC_0)$.

We consider the probability measure $q(x,\dd y)\pi_0(\dd x)$ on $\cH\times\cH$ and random variables $u\sim \pi_0$, $v = \sqrt{1-\beta^2}\,u + \beta w$ with $w\sim \cN(0,\cC_0)$ independent of $u$. Then $(u,v)$ is distributed as $q(x,\dd y)\pi_0(\dd x)$. In Example 7 of \cite[Section 5.1]{DashtiStuart:2017Bayesian} the characteristic function of $(u,v)$ is computed, and it is shown that the distribution $q(x,\dd y)\pi_0(\dd x)$ is centered Gaussian and has the same covariance operator as $q(y,\dd x)\pi_0(\dd y)$. This proves the relation \eqref{eq:auxbalance}. Finally, the special form of the acceptance probability \eqref{eq:acceptanceprob} and relation \eqref{eq:auxbalance} allow one to apply
\cite[Theorem 22]{DashtiStuart:2017Bayesian} to conclude that the condition \eqref{eq:detailedbalancepCNq} holds, hence also the desired detailed balance condition \eqref{eq:detailedbalance}. In particular, the chain produced by the pCN-algorithm has $\pi$ as invariant measure.

Convergence of the chain is a more delicate issue. It is remarked in \cite[p. 2459]{HairerStuart:2014} that the Markov chain produced by the pCN-algorithm is not $\varphi$-irreducible for whatever measure $\varphi$. Nevertheless, the authors show by different methods that an ergodic theorem holds for the pCN-algorithm under certain conditions on the potential $\Phi$ and the parameter $\beta$, at least for functions $f:\cH\to\R$ in Proposition \ref{prop:ET} that are Lipschitz with respect to a certain distance function on $\cH$ \cite[Theorem 4.3]{HairerStuart:2014}. In addition, a central limit theorem can be shown, albeit with a non-explicitly given variance $\gamma^2$
\cite[Theorem 4.4]{HairerStuart:2014}.

Other results on the Markov chain produced by the pCN-algorithm include convergence of the marginal probability distributions $p^{(n)}(x,\cdot)$ to the invariant probability measure $\pi$ in a Wasserstein metric \cite[Theorem 13]{Durmus:2016}, albeit under the condition that the potential $\Phi$ is globally H\"older continuous.
In the finite dimensional setting, various conditions on the invariant measure $\pi(\dd x) = \pi(x)\dd x$ have been highlighted in \cite{Kamatani:2017} under which the pCN-algorithm and related algorithms are geometrically ergodic or not geometrically ergodic.

\subsection{Gaussian measures in Hilbert spaces}
\label{subsec:Gaussianmeasures}

\emph{Gaussian random fields and Gaussian measures, bottom up.}  References for this material are \cite[Sections 4, 5, 7]{LordPowellShardlow:2014} and \cite[Section 37]{LoeveII:1978}. Let $D$ be a bounded domain in ${\mathbb R}^d$. A \emph{random field on $D$} is a collection of random variables
$\{\Theta(x):x\in D\}$, defined on a common probability space $(\Omega,{\mathcal A},{\rm P})$. (We assume here that every $\Theta(x)$ is real valued.) To determine the random field, one has to specify the finite dimensional joint distributions of $\Theta(x_1),\ldots, \Theta(x_n)$, $x_1, \ldots, x_n\in D$, $n\in {\mathbb N}$. The random field is \emph{Gaussian}, if all finite dimensional distributions are multivariate Gaussian.

Let $\Theta$ be a \emph{second order random field}, that is, $\Theta(x)\in L^2(\Omega)$ for all $x\in D$:
\[
  {\rm E}(|\Theta(x)|^2) = \int_\Omega |\Theta(x,\omega)|^2 {\rm P}(\dd\omega) < \infty.
\]
In this case, one defines the \emph{mean} or \emph{mean function} by $m(x) = {\rm E}(\Theta(x))$ and the \emph{covariance} or \emph{autocovariance function} by
$C(x,y) = {\rm COV}(\Theta(x),\Theta(y)) = {\rm E}((\Theta(x) - m(x))(\Theta(y)-m(y))$. The covariance is a nonnegative definite function.

A Gaussian field is completely determined by specifying the mean and the autocovariance function. This is due to the fact that the joint density of an $n$-dimensional Gaussian variable is determined by the mean vector $\mu = [\mu(x_1),\ldots, \mu(x_n)]^{\sf T}$ and the covariance matrix
\[
  {\rm C} = \begin{bmatrix}
    C(x_1,x_1) & C(x_1,x_2) & \ldots & C(x_1,x_n)\\
    \vdots &&& \vdots\\
    C(x_n,x_1) & C(x_n,x_2) & \ldots & C(x_n,x_n)
   \end{bmatrix}.
\]
In fact, the density is
\[
f(x)= \dfrac 1{(2\pi)^{n/2}\sqrt{\det {\rm C}}}\,\exp\Big(-\frac 12(x-\mu)^{\sf T}{\rm C}^{-1}(x-\mu)\Big)
\]
for ${x} = [x_1,\ldots,x_n]^{\sf T}$.

Let $\Theta$ be a random field on $D$ whose mean and autocovariance function are continuous.
Then $\Theta$ is mean-square continuous \cite[Section 37.2]{LoeveII:1978}, that is, $\Theta:D\to L^2(\Omega)$ is continuous. In particular,
\[
   \Theta\in L^2(D:L^2(\Omega)) = L^2(D\times\Omega) = L^2(\Omega:L^2(D)).
\]
Hence one may view the random field $\Theta$ as a random variable with values in the Hilbert space ${\mathcal H} = L^2(D)$.

The autocovariance function gives rise to an operator ${\mathcal C}:{\mathcal H}\to{\mathcal H}$, the \emph{covariance operator}
\[
{\mathcal C}f(x) = \int_D C(x,y)f(y){\rm d}y.
\]
A straightforward calculation shows that
\begin{equation}\label{eq:covformula}
   \langle{\mathcal C}f,g\rangle = {\rm E}(\langle\Theta-m,f\rangle\langle\Theta-m,g\rangle)
\end{equation}
for $f,g\in {\mathcal H}$.
The covariance operator is selfadjoint, positive and trace class \cite[Lemma 4.37]{LordPowellShardlow:2014}. Therefore, there is a complete orthonormal system $(\gamma_j,\varphi_j)$ on ${\mathcal H}$ with $\gamma_j \geq 0$, $\sum_{j=1}^\infty\gamma_j < \infty$, and it holds that
\[
C(x,y) = \sum_{j=1}^\infty \gamma_j\varphi_j(x)\varphi_j(y)
\]
with convergence in $L^2(D\times D)$.

Finally, one has the Karhunen-Lo\`eve expansion of the random field
\begin{equation}\label{eq:KL}
\Theta(x,\omega)  =  m(x) + \sum_{j=1}^\infty\sqrt{\gamma_j}\xi_j(\omega)\varphi_j(x)
\end{equation}
where the sum converges in $L^2(\Omega:{\mathcal H})$ and
\[
  \xi_j(\omega) = \frac{1}{\sqrt{ \gamma_j}}\int_D (\Theta(x,\omega)-m(x))\varphi_j(x){\rm d}x\quad\text{for}\ \gamma_j\neq 0;
\]
the $\xi_j, j = 1, 2, 3, \ldots$ form a sequence of orthogonal random variables (${\rm E}(\xi_j\xi_k) = \delta_{jk}$), possibly extended for $\gamma_j = 0$.

Assume now that the random field $\{\Theta(x):x\in D\}$ is Gaussian (with continuous mean and covariance function). The associated random variable $\Theta\in L^2(\Omega:{\mathcal H})$ is \emph{Gaussian} in the sense that $\langle \Theta, f\rangle = \int_D\Theta(x,\cdot)f(x){\rm d}x$ is a real valued Gaussian random variable for every
$f\in{\mathcal H}$. (Due to the continuity of $\Theta$, the integral can be approximated by Riemann sums with values in ${\mathcal H}$.) For fixed $f\in {\mathcal H}$, this Gaussian variable $\langle \Theta, f\rangle$ has mean and variance
\begin{equation}\label{eq:meanvar}
   {\rm E}(\langle \Theta, f\rangle) = \langle m, f\rangle, \qquad {\rm E}(|\langle \Theta - m, f\rangle|^2) =  \langle{\mathcal C}f,f\rangle
\end{equation}
as follows from \eqref{eq:covformula}. For later reference, we note that its characteristic function is given by
\begin{equation}\label{eq:charfun}
\chi(f) = {\rm E}\big({\rm e}^{{\rm i}\langle\Theta,f\rangle}\big) = \exp\left({\rm i}\langle m, f \rangle - \frac12\langle{\mathcal C}f,f\rangle\right).
\end{equation}
Further, the $\xi_j, j = 1, 2, 3, \ldots$ in the Karhunen-Lo\`eve expansion \eqref{eq:KL} form
a sequence of independent standard normal variables.

\emph{Gaussian random fields and Gaussian measures, top down.} Let ${\mathcal H}$ be a separable Hilbert space. We denote the $\sigma$-algebra of Borel subsets of
${\mathcal H}$ by ${\mathcal B}({\mathcal H})$ and fix a probability measure $\pi$ on $({\mathcal H},{\mathcal B}({\mathcal H}))$.

Following \cite[Section 1.4]{DaPrato:2006}, assume that $\int_{\mathcal H}|\omega|\pi(\dd\omega)$ is finite.\footnote{The Hilbert space plays the double role as a function space and as a probability space. We denote its elements by $f$ or by  $\omega$, depending on what of the two features has to be emphasized.}
Then the functional $F:{\mathcal H}\to{\mathbb R}$,
\[
   F(f) = \int_{\mathcal H}\langle\omega, f\rangle\pi(\dd\omega),\quad f\in {\mathcal H}
\]
is continuous. By the Riesz representation theorem, there exists a unique $m\in {\mathcal H}$, the \emph{mean} of $\pi$, such that
\[
   \langle m,f\rangle = \int_{\mathcal H}\langle\omega, f\rangle\pi(\dd\omega),\quad f\in {\mathcal H}.
\]
Assume further that $\int_{\mathcal H}|\omega|^2\pi(\dd\omega)$ is finite. Then the bilinear form
\[
  G(f,g) = \int_{\mathcal H}\langle\omega-m, f\rangle\langle\omega-m, g\rangle\pi(\dd\omega),\quad f, g\in {\mathcal H}
\]
is continuous. Again by the Riesz representation theorem, there is a unique linear bounded operator ${\mathcal C}: {\mathcal H}\to{\mathcal H}$ such that
\[
  \langle{\mathcal C}f,g) = \int_{\mathcal H}\langle\omega-m, f\rangle\langle\omega-m, g\rangle\pi(\dd\omega),\quad f, g\in {\mathcal H}.
\]
${\mathcal C}$ is called the \emph{covariance} or \emph{covariance operator} of $\pi$.

The covariance operator ${\mathcal C}$ is selfadjoint, positive and trace class \cite[Proposition 1.8]{DaPrato:2006}.

The probability measure $\pi$ is a \emph{Gaussian measure}, if its characteristic function $\chi(f) = \int_{\mathcal H}{\rm e}^{{\rm i}\langle\omega, f\rangle}\pi(\dd\omega)$
is of the form
\begin{equation}\label{eq:charfunGauss}
\chi(f) = \exp\left({\rm i}\langle m, f \rangle - \frac12\langle{\mathcal C}f,f\rangle\right),\quad f\in {\mathcal H}.
\end{equation}
With somewhat more effort one can show that, given $m\in{\mathcal H}$ and a selfadjoint, positive trace class operator ${\mathcal C}$ on ${\mathcal H}$, there is a unique Gaussian measure $\pi$ on $({\mathcal H},{\mathcal B}({\mathcal H}))$ with characteristic function \eqref{eq:charfunGauss} \cite[Theorem 1.12]{DaPrato:2006}.

Let $(\Omega, {\mathcal A},{\rm P})$ be a probability space and $\Theta: \Omega\to {\mathcal H}$ be an ${\mathcal H}$-valued random variable, i.e.,
$({\mathcal A},{\mathcal B}({\mathcal H}))$-measurable. Recall that the \emph{law of $\Theta$} is the measure $\pi = \Theta_\sharp{\rm P}$ on ${\mathcal H}$, given by $\pi(B) = {\rm P}(\Theta^{-1}(B))$ for $B \in{\mathcal B}({\mathcal H}) $. The random variable $\Theta$ is called \emph{Gaussian}, if its law is a Gaussian measure on ${\mathcal H}$.
We denote this by $\Theta \sim{\mathcal N}(m,{\mathcal C})$.

{\emph{Examples.} (a) Any Gaussian measure on a separable Hilbert space ${\mathcal H}$ gives rise to a Gaussian random variable $\Theta_0: \Omega = {\mathcal H} \to {\mathcal H}$, given as the identity map $\Theta_0(\omega) = \omega$.

(b) The random variable $\Theta$ on ${\mathcal H} = L^2(D)$ corresponding to a Gaussian random field $\{\Theta(x): x\in D\}$ is Gaussian, thanks to \eqref{eq:charfun}. Indeed, let $\cC$ be the operator with integral kernel $C(x,y)$. We have the formulas \eqref{eq:meanvar} for the mean and variance for
$\langle \Theta, f \rangle$. The characteristic function for $\Theta_\sharp{\rm P}$ is
\[
   \int_{\mathcal H}{\ee}^{{\rm i}\langle\omega, f\rangle}\Theta_\sharp{\rm P}(\dd\omega)
      = \int_{\mathcal H}{\ee}^{{\rm i}\langle\Theta(\omega), f\rangle}P(\dd\omega)
      = {\rm E}\big({\rm e}^{{\rm i}\langle\Theta,f\rangle}\big) = \exp\left({\rm i}\langle m, f \rangle - \frac12\langle{\mathcal C}f,f\rangle\right)
\]
by \eqref{eq:charfun}. Thus the law of $\Theta$ is the Gaussian measure $\pi$ with mean $m$ and covariance operator $\cC$.

Let us revisit the Karhunen-Lo\`eve expansion in the setting of a Gaussian measure $\pi$ on a Hilbert space $\cH$.
The covariance operator $\cC$ has a complete orthonormal system $(\gamma_j,\varphi_j)$ on ${\mathcal H}$ of eigenvalues $\gamma_j$ and eigenfunctions $\varphi_j$ with $\gamma_j \geq 0$, $\sum_{j=1}^\infty\gamma_j < \infty$.
Let $(\Omega,\cA,P)$ be a probability space with a sequence $\xi_j$, $j = 1,2,3,\ldots$ of independent identically distributed standard Gaussian variables. Set
\[
   \Theta_N(\omega) =  m + \sum_{j=1}^N\sqrt{\gamma_j}\xi_j(\omega)\varphi_j, \quad\langle\Theta_N(\omega),f\rangle = \langle m,f\rangle + \sum_{j=1}^N\sqrt{\gamma_j}\xi_j(\omega)\langle\varphi_j,f\rangle
\]
for $f\in\cH$. It is clear that $\Theta_N: \Omega \to \cH$ is measurable and the $\langle\Theta_N(\omega),f\rangle$ are real-valued Gaussian variables with mean $\langle m,f\rangle$ and variance
\[
   \rE\Big(\sum_{j=1}^N\sum_{k=1}^N\sqrt{\gamma_j}\xi_j(\omega)\langle\varphi_j,f\rangle\sqrt{\gamma_k}\xi_k(\omega)\langle\varphi_k,f\rangle\Big)
   = \sum_{j=1}^N{\gamma_j}\langle\varphi_j,f\rangle^2 \leq \sum_{j=1}^N{\gamma_j} \|f\|_\cH^2.
\]
The convergence of the variances of a sequence of Gaussian variables implies its mean-square convergence. Thus the $\langle\Theta_N(\omega),f\rangle$ converge
in $L^2(\Omega)$ to a Gaussian random variable. The convergence is uniform on the unit ball of $\cH$, so the $\Theta_N$ converge to a random variable $\Theta\in L^2(\Omega:\cH)$.
Clearly,
\[
   \rE(\langle \Theta, f\rangle) = \langle m, f \rangle,\quad {\rm E}(|\langle \Theta - m, f\rangle|^2)
     = \sum_{j=1}^\infty{\gamma_j}\langle\varphi_j,f\rangle^2 =  \langle{\mathcal C}f,f\rangle.
\]
As above, this leads to the characteristic function \eqref{eq:charfunGauss} and shows that the law of $\Theta$ is the Gaussian measure $\pi$.

Summarizing, given a Gaussian measure $\pi$ on a Hilbert space $\cH$, the Karhunen-Lo\`eve expansion yields a random variable $\Theta$ on $\cH$ whose law is $\pi$. This opens the possibility of producing samples of $\Theta$ by drawing from a sequence of independent standard Gaussian random variables.

Suppose $\cH = L^2(D)$ for some domain $D\subset\R^d$ and the realizations of the random variable $\Theta\in L^2(\Omega:\cH)$ are almost surely continuous functions on $D$. Then $\Theta$ corresponds to a classical Gaussian random field on $D$, which has $\cC$ as its covariance operator.
In this sense, the bottom-up approach and the top-down approach are equivalent. The bottom-up approach offers the possibility to start with some explicitly given autocovariance function, but the explicit calculation of the orthonormal system might be difficult, see the examples in \cite[Section 2.3.3]{Ghanem:2003}.
The top-down approach offers the possibility to start from a simple definition of the covariance operator, e.g. as a negative fractional power of an elliptic partial differential operator $\cL$  (often referred to as the \emph{precision operator}), whose eigenfunction expansion is well-known. Then one has less control over the kernel $C(x,y)$, but one may work directly with the Karhunen-Lo\`eve expansion without computing the kernel.

\emph{Example 1.} To illustrate the top-down approach, consider the Hilbert space ${\mathcal H} = L^2(D)$ on the domain $D = [0,1]^d\subset {\mathbb R}^d$. The operator ${\mathcal L}$ on ${\mathcal H}$ is the Laplacian $-\Delta$ with operator domain ${\mathcal D}({\mathcal L}) = H^2(D)\cap H_0^1(D)$. The operator ${\mathcal L}$ is selfadjoint, positive and an isomorphism between ${\mathcal D}({\mathcal L})$ and ${\mathcal H}$. There is a complete orthonormal system $(\lambda_k, \varphi_k)$, indexed by $k = (k_1,\ldots, k_d)\in{\mathbb N}^d$ with eigenvalues and eigenfunctions given by
\[
    \lambda_k = \pi^2 |k|^2 = \pi^2(k_1^2 + \cdots + k_d^2),\qquad \varphi_k(x) = \prod_{i=1}^d\big(\sqrt{2}\sin(k_i\pi x_i)\big).
\]
We have that ${\mathcal L}g = \sum_{k\in{\mathbb N}^d} \lambda_k\langle g,\varphi_k\rangle\varphi_k$. For $\alpha \geq 0$, the sequence $\lambda_k^{-\alpha}$ is bounded and so
${\mathcal C} = {\mathcal L}^{-\alpha}$ maps ${\mathcal H} = L^2(D)$ continuously into itself via
${\mathcal C}f = \sum_{k\in{\mathbb N}^d} \lambda_k^{-\alpha}\langle f,\varphi_k\rangle\varphi_k$. The operator ${\mathcal C}$ is selfadjoint and positive. Further, for $\alpha > d/2$ it is trace class, because $\sum_{k\in{\mathbb N}^d} \lambda_k^{-\alpha} = \pi^{-2\alpha}\sum_{k\in{\mathbb N}^d}|k|^{-2\alpha}$ is finite. Thus ${\mathcal C}$ is the covariance operator of a Gaussian measure on ${\mathcal H}$. The operator ${\mathcal C}$ is an integral operator with kernel
\[
   C(x,y) = \sum_{k\in{\mathbb N}^d} \lambda_k^{-\alpha} \varphi_k(x)\varphi_k(y) \in L^2(D\times D).
\]
The eigenfunctions $\varphi_k$ are uniformly bounded by $2^{d/2}$ and the sequence $\lambda_k^{-\alpha}$ is summable, thus the convergence is actually uniform in $(x,y)$ and the kernel is a continuous function on $D\times D$.

On the other hand, let $(\Omega,{\mathcal A},{\rm P})$ be a probability space on which a sequence $(\xi_k, k \in {\mathbb N}^d)$ of independent standard Gaussian random variables is defined. We can introduce a random field on $D$ by
\[
\Theta(x,\omega) = \sum_{k\in{\mathbb N}^d}\lambda_k^{-\alpha/2}\xi_k(\omega)\varphi_k(x).
\]
We have that
\[
   {\rm E}(\Theta(x)^2)
   = {\rm E}\big(\sum_{k\in{\mathbb N}^d}\lambda_k^{-\alpha/2}\xi_k(\omega)\varphi_k(x)\sum_{l\in{\mathbb N}^d}\lambda_l^{-\alpha/2}\xi_l(\omega)\varphi_l(x)\big)
   = \sum_{k\in{\mathbb N}^d}\lambda_k^{-\alpha}\varphi_k(x)^2.
\]
By the considerations above, the sequence converges uniformly on $D$, thus ${\rm E}(\Theta(x)^2)$ is a continuous function on $D$. Clearly, the covariance operator generated by $\Theta$ is ${\mathcal C}$. For $x\in \partial D$,
${\rm E}(\Theta(x)^2) = 0$ and the values of any finite partial sum on the boundary $\partial D$ are zero almost surely.
When applied in the reconstruction of random fields in Bayesian inverse problems, this is undesirable.

Accordingly, we will change the setting to the Laplace operator with periodic boundary conditions next.

\emph{Example 2.} Consider the Hilbert space ${\mathcal H} = L^2(0,1)$. The operator ${\mathcal L}$ will be the Laplacian $-\partial_x^2$ with operator domain
${\mathcal D}({\mathcal L}) = H^2_{\rm per}(0,1) = \{f\in H^2(0,1): f(0) = f(1), f'(0) = f'(1)\}$. The operator ${\mathcal L}$ is selfadjoint, positive and an isomorphism between ${\mathcal D}({\mathcal L})$ and its range ${\mathcal R}({\mathcal L}) = \{g\in L^2(D): \int_0^1 g(x){\rm d}x = 0\}$.
The operator ${\mathcal L}$ has a sequence of eigenvalues $\lambda_0 = 0$, $\lambda_k = (2\pi k)^2$ of multiplicity two for $k\geq 1$ with eigenfunctions $\psi_0 = 1$,
$\varphi_k(x) = \sqrt{2}\sin 2k\pi x$, $\psi_k(x) = \sqrt{2}\cos 2k\pi x$.

The fact that ${\mathcal H} = ({\mathbb R}\cdot 1) \oplus {\mathcal R}({\mathcal L})$ suggests to extend ${\mathcal L}^{-1}$ to $(c + {\mathcal L}^{-1})$,
\begin{eqnarray*}
 (c + {\mathcal L}^{-1})f(x) &=&
     c\langle\psi_0, f \rangle + \sum_{k=1}^\infty (2k\pi)^{-2}\Big (\langle \varphi_k, f \rangle \varphi_k(x) + \langle \psi_k, f \rangle\psi_k(x)\Big)
\end{eqnarray*}
which is a positive and selfadjoint operator if $c>0$. For $\alpha > 1/2$, the operator
\begin{eqnarray}\label{eq:kernel}
{\mathcal C}f(x) \ =\  (c + {\mathcal L}^{-1})^\alpha f(x) &=&
     c^\alpha\langle\psi_0, f \rangle + \sum_{k=1}^\infty (2k\pi)^{-2\alpha}\Big (\langle \varphi_k, f \rangle \varphi_k(x) + \langle \psi_k, f \rangle\psi_k(x)\Big)\nonumber\\
   &=& \int_0^1 \Big( c^\alpha + 2 \sum _{k=1}^\infty (2k\pi)^{-2\alpha}\cos 2 k\pi (x-y)\Big)   f(y){\rm d}y
\end{eqnarray}
is trace class and can serve as the covariance operator of a Gaussian measure on ${\mathcal H}$. The last equality follows by expanding the terms and allows one to read off the integral kernel.

Let $(\Omega,{\mathcal A},{\rm P})$ be a probability space with a double sequence $\xi_k, k\geq 1$, $\eta_k, k \geq 0$ of independent standard Gaussian random variables.
The random field
\[
\Theta(x,\omega) = c^{\alpha/2}\eta_0(\omega) + \sqrt{2}\sum_{k=1}^\infty (2k\pi)^{-\alpha}\big( \xi_k(\omega)\sin 2k\pi x + \eta_k(\omega)\cos 2k\pi x\big)
\]
is Gaussian and has ${\mathcal C}$ as covariance operator.

\emph{Remark.} For $\alpha = 1$ one can get an explicit formula for the kernel $C(x,y)$ of ${\mathcal C}$. Comparing with the Fourier series for
\[
    \frac12\Big(z-\frac12\Big)^2 = \frac{1}{24} + 2 \sum _{k=1}^\infty (2k\pi)^{-2}\cos 2 k\pi z
\]
we find that
\begin{eqnarray*}
  C(x,y) = c + 2 \sum _{k=1}^\infty (2k\pi)^{-2}\cos 2 k\pi (x-y) = \frac12\Big(x-y-\frac12\Big)^2 + c - \frac{1}{24}.
\end{eqnarray*}
In the bottom-up approach, we would have started with the Fourier expansion of $C(x,0)$ and then would have worked upwards to construct a Gaussian, trigonometric random field.

\emph{Example 3.} In the application to reconstructing the right-hand side in the two-dimensional Poisson equation in Subsection \ref{subsec:2dPoisson}, we extend Example 2 to two dimensions, using the Gaussian random field on $D = [0,1]\times[0,1]$,
\[
   \Theta(x_1,x_2,\omega) =
\xi_0(\omega) + 2(2\pi)^{-\alpha}\!\!\!\sum_{\substack{k,l = 0\\(k,l)\neq(0,0)}}^\infty\sum_{i=1}^4 (k^2 + l^2)^{-\alpha/2}\xi_{kl}^{(i)}\psi_{kl}^{(i)}(x_1,x_2)
\]
with $\alpha > 1$, where
\begin{eqnarray*}
\psi_{kl}^{(1)}(x_1,x_2) = \sin 2k\pi x_1\sin 2l\pi x_2, \quad\psi_{kl}^{(2)}(x_1,x_2) = \sin 2k\pi x_1\cos 2l\pi x_2, \\
\psi_{kl}^{(3)}(x_1,x_2) = \cos 2k\pi x_1\sin 2l\pi x_2, \quad\psi_{kl}^{(4)}(x_1,x_2) = \cos 2k\pi x_1\cos 2l\pi x_2
\end{eqnarray*}
and the $\xi_0$ as well as  $\xi_{kl}^{(i)}$, $k,l\in{\mathbb N}, i = 1,\ldots, 4$ are all independent standard Gaussian random variables.

For reference in Subsection \ref{subsec:2dPoisson} we denote the corresponding covariance operator by $\cC_\alpha$.

\emph{Example 4.} The covariance operator may be taken as a negative power of a general second order elliptic partial differential operator $\cL$ satisfying suitable conditions
\cite[Assumption 2.9]{Stuart_2010_acta_numerica}. This gives us the opportunity to have a glimpse at a different road to covariances, the \emph{SPDE approach} \cite{Lindgren:2022}. Informally speaking, if a random field $\Theta$ on a domain $D$ is defined as $\Theta = \cL^{-1}\dot{W}$, where $\dot{W}$ is Gaussian white noise on $D$, then $\Theta$ has covariance operator $\cL^{-2}$. To exemplify this, consider the Sturm-Liouville operator
\[
   \cL f(x) = - \frac{\dd}{\dd x}\left(p(x)\frac{\dd}{\dd x}f(x)\right) + q(x)f(x)
\]
on $D = (0,1)$ where $p(x) > 0$, $q(x)\geq 0$. We assume that the coefficient functions $p,q$ are sufficiently regular so that $\cL$ is an isomorphism from
$H^1_0(D)$ to $H^{-1}(D)$. Gaussian white noise $\dot{W}$ on $D = (0,1)$ can be seen as the pathwise weak derivative of Wiener process, thus certainly belongs to $H^{-1}(D)$. The equation $\cL(\Theta) = \dot{W}$ has a unique solution $\Theta \in H^1_0(D)$.
It can be written as an It\^o integral in terms of the Green function $G$ of $\cL$,
\[
   \Theta(x) = \int_0^1 G(x,z)\dd W(z).
\]
The random field $\Theta$ has mean zero and covariance (using the It\^o isometry)
\[
  \rE(\Theta(x)\Theta(y)) = \rE\left(\int_0^1 G(x,z)\dd W(z)\int_0^1 G(y,z)\dd W(z)\right) = \int_0^1 G(x,z)G(y,z)\dd z.
\]
By the symmetry of Green's function, the latter term equals $\int_0^1 G(x,z)G(z,y)\dd z$ and hence coincides with the integral kernel of $\cL^{-2}$.

The law of $\Theta$ is $\cN(0,\cL^{-2})$. Sampling from one-dimensional white noise or from the solution of an It\^o equation is easy, thus the SPDE approach is a viable alternative to the Karhunen-Lo\`{e}ve expansion. A further advantage in higher space dimensions is that the shape of the domain $D\subset\R^d$ may be arbitrary; no eigenvalues need to be computed. However, regularity questions enter. Spatial white noise is a generalized stochastic process, defined through its action on test functions $\varphi$ (smooth functions of compact support in $D$), $\varphi\to \dot{W}(\omega)(\varphi)$, $\omega\in\Omega$. It is well-known (\cite[Section 2.1]{HoldenOksendal:2010}) that the action of $\dot{W}$ can be extended to $\varphi\in L^2(D)$.
Then the Green's function argument goes through as long as $G(x,\cdot) \in L^2(D)$ for every $x\in D$. If $\cL = -\Delta$ is the Laplacian, this is the case for space dimensions $d = 1,2,3$; see also the discussion in \cite{Stadler:2013}.

\emph{Remarks.} (a) The concept of the SPDE approach goes back at least to \cite{Whittle:1954} and was notably used to construct random fields with a Mat\'{e}rn covariance $C(x,y) \sim (\kappa(|x-y|)^\nu K_\nu(\kappa|x-y|)$ as solutions to $(\kappa^2 - \Delta)^{\alpha/2}\Theta = \dot{W}$ on $\R^d$, where $\kappa, \nu > 0$, $\alpha = \nu + d/2$ and $K_\nu$ is the modified Bessel function of the second kind. We refer to \cite{Lindgren:2011} and the survey paper \cite{Lindgren:2022}, which also discusses an extension to the non-constant, anisotropic case obtained through $$\tau(x)\big(\kappa(x)^2 - \nabla\cdot H(x)\nabla\big)^{\alpha/2}\Theta(x) = \dot{W}(x).$$ Further, $\kappa(x)$ may be taken as a random field itself, as done in \cite{Roininen:2019} with $\alpha = 2$. Mat\'{e}rn covariances together with a Karhunen-Lo\`{e}ve expansion were also used in \cite{Taghizadeh:2025}.

(b) Sampling from a random field on arbitrary domains has its difficulties, due to boundary effects. This can be mitigated either by extending the domain to a box or by optimizing the boundary conditions; this issue is discussed and investigated in \cite{Daon:2018}.

(c) Apart from papers quoted in (a), non-constant coefficient precision operators $\cL$ arise in inverse problems on manifolds \cite{Harlim:2022,Lindgren:2011,Petra:2014}, e. g., as Laplace-Beltrami operators.

(d) A further method for drawing samples from a Gaussian random field is the use of so-called \emph{Langevin algorithms}. These algorithms consist in constructing Langevin-type diffusion equations having the desired measure as stationary distribution \cite[Section 7.8.5]{robert1999monte}, \cite{HairerStuartVoss:2007,HairerStuartVoss:2009,HairerStuartVoss:2005}.

\end{document}